\documentclass[11pt,letterpaper]{article}

\usepackage[T1]{fontenc}
\usepackage[utf8]{inputenc}
\usepackage{lmodern}
\usepackage{amsmath,amssymb,amsthm,mathtools}
\usepackage{microtype}
\usepackage{xcolor}
\usepackage[margin=1in]{geometry}
\usepackage[colorlinks=true,linkcolor=blue!55!black,citecolor=blue!55!black,urlcolor=blue!55!black]{hyperref}
\hypersetup{
  pdftitle={Sharp Last-Iterate Bounds for the Stochastic subGradient Method in Fixed Dimension},
  pdfsubject={Last-iterate error under conditional noise moments in finite dimension}
}
\usepackage{authblk}

\allowdisplaybreaks
\newtheorem{theorem}{Theorem}
\newtheorem{lemma}[theorem]{Lemma}

\newcommand{\R}{\mathbb R}
\newcommand{\E}{\mathbb E}
\newcommand{\Pp}{\mathbb P}
\newcommand{\cC}{\mathcal C}
\newcommand{\cF}{\mathcal F}
\newcommand{\cG}{\mathcal G}

\newcommand{\cI}{\mathcal I}
\newcommand{\cX}{\mathcal X}
\newcommand{\1}{\mathbf 1}
\newcommand{\dist}{\operatorname{dist}}

\DeclareMathOperator*{\argmin}{argmin}

\title{Sharp Last-Iterate Bounds for the Stochastic subGradient Method}

\author[1]{Guglielmo Beretta}
\author[2]{Tommaso Cesari}
\author[3]{Roberto Colomboni}
\author[1]{Andrea Paudice}

\affil[1]{Department of Computer Science, Aarhus University, Aarhus, Denmark}
\affil[2]{School of Electrical Engineering and Computer Science, University of Ottawa, Ottawa, Canada}
\affil[3]{School of Mathematics, University of Bristol, Bristol, United Kingdom}

\affil[ ]{
\texttt{gberetta@cs.au.dk},
\texttt{tcesari@uottawa.ca},
\texttt{roberto.colomboni@bristol.ac.uk},
\texttt{apaudice@cs.au.dk}
}

\begin{document}
\maketitle

\begin{abstract}
We determine the worst-case expected optimization error of the last iterate of constant-stepsize projected Stochastic
subGradient Method on convex Lipschitz objectives. The upper and lower bounds match up to universal
constants for every stepsize, horizon, dimension, and conditional noise moment of order $p\in[1,\infty]$.
With bounded noise and a fixed positive bound on the initial distance to a minimizer, the best constant
stepsize gives worst-case error after $n$ iterates of order $1/\sqrt n$ on the line and
$\sqrt{\log n/(n\log\log n)}$ in every fixed dimension at least two.
\end{abstract}

\section{Problem formulation}
\label{sec:model}

Let the dimension $d\ge1$ and the horizon $n\ge2$ be integers. On $\R^d$ we consider the Euclidean inner
product $\langle\cdot,\cdot\rangle$ and the Euclidean norm $\|\cdot\|$. Let the feasible set $\cX$ be a nonempty
closed convex subset of $\R^d$. Fix $L\ge0$ and let $f:\cX\to\R$ be convex and $L$-Lipschitz, so that for all
$x,y\in\cX$ we have $|f(x)-f(y)|\le L\|x-y\|$. Assume that $f$ attains its minimum
$f_\star\coloneqq\min_{x\in\cX}f(x)$, let $\cX_\star\coloneqq\argmin_{x\in\cX}f(x)$ be the set
of minimizers, and for each $x\in\cX$ let $\Delta(x)\coloneqq f(x)-f_\star$ be its associated optimization error.

Let $\Pi_{\cX}$ denote Euclidean projection onto $\cX$. For $x\in\cX$, a vector $g\in\R^d$ belongs
to the subdifferential $\partial_{\cX}f(x)$ if, for every $y\in\cX$, we have
$f(y)\ge f(x)+\langle g,y-x\rangle$. When $\cX=\R^d$ we write $\partial f$. Let $X_1\in\cX$ be
deterministic and let $D_1\coloneqq\min_{x\in\cX_\star}\|X_1-x\|$ be its distance to the minimizers.

Fix a constant stepsize $\eta>0$, a noise moment order $p\in[1,\infty]$, and a noise scale $\sigma\ge0$.
On a probability space with probability $\Pp$ and expectation $\E$, let $(\cF_k)_{k=1}^n$ be a
filtration. For each $1\le k<n$, given the $\cF_k$-measurable iterate $X_k$, choose an
$\cF_k$-measurable subgradient $G_k\in\partial_{\cX}f(X_k)$ with $\|G_k\|\le L$.
Let the noise $W_k\in\R^d$ be $\cF_{k+1}$-measurable and satisfy, almost surely,
\begin{equation}\label{eq:conditional-moment}
 \E[W_k\mid\cF_k]=0,\qquad
 \begin{cases}
  \E[\|W_k\|^p\mid\cF_k]\le\sigma^p,&\textrm{if }p<\infty,\\
  \|W_k\|\le\sigma,&\textrm{if } p=\infty .
 \end{cases}
\end{equation}
The noise law may therefore depend on the history. We define the next iterate by
\begin{equation}\label{eq:ssgm}
 X_{k+1}\coloneqq\Pi_{\cX}\bigl(X_k-\eta(G_k+W_k)\bigr).
\end{equation}

For $D\ge0$, let $\cC_p^{(d)}(D)$ consist of the instances above with $L,\sigma\le1$ and $D_1\le D$.
The feasible set may be unbounded. When $D=0$, the method starts at a minimizer.
An instance $\cI$ may depend on $n$ and $\eta$. Writing $\E_{\cI}$ for its expectation, define the worst-case expected last-iterate error by
\begin{equation}\label{gd:eq:risk}
 \mathcal R_{n,p}^{(d)}(\eta,D)\coloneqq
 \sup_{\cI\in\cC_p^{(d)}(D)}\E_{\cI}[\Delta(X_n)] .
\end{equation}

\section{Main theorem}
\label{sec:main}
Using natural logarithms and writing $\mathbf1$ for an indicator, let
\begin{equation}\label{eq:factor}
 b_{n,p,d}\coloneqq\max\!\left\{
  1+\frac{n^{2-p}-1}{2-p},\quad
  \mathbf 1_{\{d\ge2\}}\,\frac{\log n}{\log(e+(\log n)/d)}
 \right\},
\end{equation}
where, at $p=2$ and $p=\infty$, we define the first term in the maximum by continuity as $1+\log n$ and $1$, respectively. We use $\lesssim$, $\gtrsim$, and
$\asymp$ for comparisons up to positive universal constants, with a subscript for any dependence on fixed parameters.

\begin{theorem}\label{thm:main}
For every integer $d\ge1$, every integer $n\ge2$, every $p\in[1,\infty]$, every $D\ge0$, and every
$\eta>0$, the worst-case expected last-iterate error satisfies
\begin{equation}\label{eq:main-prescribed}
 \mathcal R_{n,p}^{(d)}(\eta,D)
 \asymp \min\!\left\{D,\frac{D^2}{\eta n}\right\}+\eta b_{n,p,d}.
\end{equation}
For $D>0$, the stepsize $\eta=D/\sqrt{n b_{n,p,d}}$ attains the optimal order
\begin{equation}\label{eq:main-optimized}
 \inf_{\eta>0}\mathcal R_{n,p}^{(d)}(\eta,D)
 \asymp D\sqrt{\frac{b_{n,p,d}}{n}} .
\end{equation}
\end{theorem}

For fixed $p$ and $d$ and positive $D$, substituting \eqref{eq:factor} into \eqref{eq:main-optimized} gives, as $n\to\infty$,
\[
 \inf_{\eta>0}\mathcal R_{n,p}^{(d)}(\eta,D)\asymp_{p,d}D
 \begin{cases}
  1,&p=1,\\
  n^{-(p-1)/2},&1<p<2,\\
  \sqrt{(\log n)/n},&p=2,\\
  n^{-1/2},&2<p\le\infty,\ d=1,\\
  \sqrt{\log n/(n\log\log n)},&2<p\le\infty,\ d\ge2 .
 \end{cases}
\]
With bounded noise and $D=1$, the optimized planar rate in the last line also holds on the Euclidean
unit ball (Appendix~\ref{gd:lb:localization}).
When $p$ or $d$ varies with $n$, \eqref{eq:main-optimized} retains universal constants throughout.
For example, in the plane the transition to the bounded-noise rate occurs at
$p-2\asymp(\log\log n)/\log n$, as shown in Appendix~\ref{app:transition}.

\section*{AI Disclosure}

We used OpenAI's GPT-6 Astra to brainstorm and refine algorithmic ideas, develop preliminary proof arguments, and assist with \LaTeX{} formatting.
Since the authors are not native English speakers, we also used these models to improve grammar, punctuation, and phrasing throughout the manuscript.

\clearpage
\appendix
\numberwithin{theorem}{section}
\numberwithin{equation}{section}

\section{Proof of Theorem~\ref{thm:main}}
\label{app:main}

\begin{proof}
We first keep $\eta$ fixed and account for three ways in which the terminal error can be large.
A small step can leave the starting point far from a minimizer, a late noise increment can undo
progress on a line, and in dimension at least two bounded noise can move the iterate through regions
of successively larger objective value. The lower bound will use a different instance for each effect.

To refer to each term of \eqref{eq:factor} separately in the proof, let $\Phi_{n,p}\coloneqq1+\int_1^n t^{1-p}\,dt$ with
$\Phi_{n,\infty}\coloneqq1$. For $d\ge2$ let $q_{n,d}\coloneqq(\log n)/\log(e+(\log n)/d)$.
These are the two quantities in \eqref{eq:factor}. Let
$a_\eta\coloneqq\min\{D,D^2/(\eta n)\}$ denote the contribution of the initial distance.

We begin with the upper bound. In every dimension, Appendix~\ref{gd:app:upper} shows that the
conditional variance at update $k$ enters a last-iterate bound with weight $1/(n-k)$. When $p\le2$,
truncating the noise at a level proportional to $n-k$ and controlling the change in the trajectory
gives $\mathcal R_{n,p}^{(d)}(\eta,D)\le a_\eta+8\eta\Phi_{n,p}$ in \eqref{gd:eq:heavy-fixed-step}.
In this range $\Phi_{n,p}\ge1+\log n$, so the moment term already controls $q_{n,d}$.
On the line, the two endpoints of a sublevel interval determine the directions in which the objective
can increase. Following excursions past these endpoints, Appendix~\ref{sec:upper} gives the bound
$a_\eta+(65/2)\eta\Phi_{n,p}$ for every $p$ in \eqref{eq:normalized-upper}.
For $d\ge2$ and $p>2$, Appendix~\ref{app:geometric} controls the increase over nested time intervals
ending at $X_n$. It charges the large noise increments to $\Phi_{n,p}$ and bounds the remaining
increase by $q_{n,d}$, which gives the upper bound in \eqref{eq:main-prescribed}.

We next show that each contribution is necessary. In Appendix~\ref{sec:lower}, a noiseless linear
objective gives error at least $a_\eta/4$, and an independent noise sequence on a half-line gives
error at least $\eta\Phi_{n,p}/128$ from a minimizer. Appending zero coordinates preserves both
examples. To obtain the dimension term, Appendix~\ref{gd:lb:appendix} constructs a planar trajectory
that gains a fixed amount of objective value at each successful stage. Appendix~\ref{gd:app:joint}
runs this construction on successive coordinate pairs and retains the gains from earlier pairs.
For all sufficiently large $n$, \eqref{gd:eq:joint-lower} gives a lower bound of order $\eta q_{n,d}$
with a universal constant. For the remaining horizons $q_{n,d}\le\log n$ is bounded by a numerical
constant, so the scalar example suffices. Both stochastic examples start at a minimizer, and the
bounded-noise example satisfies every moment assumption. Since the risk is a supremum and the sum of
two nonnegative numbers is at most twice their maximum, these constructions give the lower bound in
\eqref{eq:main-prescribed}, including when $D=0$.

It remains to choose the stepsize. Since $1\le b_{n,p,d}\le n$, for $\eta n\le D$ we have
$a_\eta=D\ge D\sqrt{b_{n,p,d}/n}$. For $\eta n>D$, the geometric mean of $D^2/(\eta n)$ and
$\eta b_{n,p,d}$ is $D\sqrt{b_{n,p,d}/n}$, so \eqref{eq:main-prescribed} gives the same lower
order for every stepsize. At $\eta=D/\sqrt{n b_{n,p,d}}$ the two terms are equal, which proves
\eqref{eq:main-optimized}.
\end{proof}

\subsection{The transition above the second moment}
\label{app:transition}

We now evaluate \eqref{eq:factor} when the available moment approaches two with the horizon.
For $p=2+\varepsilon$ with $0\le\varepsilon\le1$, let $u\coloneqq\log n$. Since for $x\ge0$,
$(1-e^{-1})\min\{x,1\}\le1-e^{-x}\le\min\{x,1\}$, we have
\[
 \Phi_{n,2+\varepsilon}
 =1+\frac{1-e^{-\varepsilon u}}{\varepsilon}
 \asymp \frac{1+u}{1+\varepsilon u},
\]
where the quotient in the middle is interpreted by continuity at zero. In the plane the dimension
term is $u/\log(e+u/2)$. Thus when $\varepsilon u$ stays bounded, $b_{n,p,2}$ has order $u$ as
$n\to\infty$. When $1\ll\varepsilon u\ll\log u$, the moment term has order $1/\varepsilon$ and
dominates. Once $\varepsilon u\gtrsim\log u$, the dimension term gives order $u/\log u$.
In any dimension $d\ge\log n$ with $d\ge2$, the denominator in $q_{n,d}$ is bounded and the
stochastic factor has order $1+\log n$ for all $p\ge2$.

There is a corresponding sharp answer if one suppresses the moment term and allows a constant
depending on $p$. For $2<p<\infty$, let $K_p$ be the least constant such that for every $d\ge2$,
$n\ge2$, and $\eta>0$, $\mathcal R_{n,p}^{(d)}(\eta,0)\le K_p\eta q_{n,d}$. Since $q_{n,d}$
increases with $d$, \eqref{eq:main-prescribed} reduces the calculation to the supremum of
$\max\{1,\Phi_{n,p}/q_{n,2}\}$ over $n\ge2$. Let $a\coloneqq1/(p-2)$.
When $a\ge1$, we have $\Phi_{n,p}\le1+\min\{u,a\}$. For $u\le a$, dividing by $q_{n,2}$ gives
an upper bound of order $\log(e+a)$, and for $u>a$ the same estimate follows since
$\log(e+u/2)/u$ decreases with $u$. Choosing $n=\lceil e^a\rceil$ gives $u\asymp a$ and
$\Phi_{n,p}\ge(1-e^{-1})a$, which gives the reverse inequality. When $a<1$, $\Phi_{n,p}\le2$
and $q_{n,2}$ is bounded below by a positive numerical constant. Together these cases give
$K_p\asymp\log(e+1/(p-2))$.

\section{Upper bounds in arbitrary dimension}
\label{gd:app:upper}

Our aim is to bound the last error even when the noise has no second moment. We first determine how much
the variance at each update contributes when all variances are finite. We will then modify large noise
increments, using those weights to choose how much noise to keep at each update. We take empty sums to be zero and empty products to be one. For positive integers $j$ and
finite $p\ge1$, let
\[
 H_j\coloneqq\sum_{i=1}^j\frac1i,\qquad
 S_{j,p}\coloneqq\sum_{i=1}^j i^{1-p}.
\]
Since $t\mapsto t^{1-p}$ is nonincreasing on $[1,\infty)$, the sum $S_{n-1,p}$ lies between
$\max\{1,\int_1^n t^{1-p}\,dt\}$ and $1+\int_1^{n-1} t^{1-p}\,dt$. In particular,
\begin{equation}
 \tfrac12\Phi_{n,p}\le S_{n-1,p}\le\Phi_{n,p}\le n .
 \label{eq:sum-comparison}
\end{equation}

\subsection{A bounded measurable subgradient at every feasible point}

The model specifies the subgradients along the original trajectory. To run an auxiliary process on the
same objective, we also need a bounded selection at points that the original trajectory may never visit.
We construct a Borel map $v:\cX\to\R^d$ such that, for every $x\in\cX$, $v(x)\in\partial_{\cX}f(x)$ and
$\|v(x)\|\le L$.

At a boundary point $x$ of $\cX$, the subdifferential $\partial_{\cX}f(x)$ can contain vectors of arbitrarily large norm. To
obtain a bounded choice everywhere, we extend $f$ to a convex function on $\R^d$ with the same Lipschitz
bound. For $x\in\R^d$, let $\bar f(x)\coloneqq\inf_{y\in\cX}\{f(y)+L\|x-y\|\}$. Fix $y_0\in\cX$. By the
Lipschitz and triangle inequalities, every term in the infimum is at least $f(y_0)-L\|x-y_0\|$, while
choosing $y=y_0$ gives a finite upper bound. Taking convex combinations of competing points in $\cX$
proves convexity, and the triangle inequality proves $L$-Lipschitz continuity. For $x\in\cX$, Lipschitz continuity gives $f(y)+L\|x-y\|\ge f(x)$, with equality at $y=x$, so $\bar f$
agrees with $f$ there. At every $x\in\R^d$, the finite convex
function $\bar f$ has a nonempty subdifferential $\partial\bar f(x)$, and its Lipschitz bound gives norm
at most $L$ for every subgradient.

For $x\in\R^d$, let $v(x)$ be the unique minimum-norm element of $\partial\bar f(x)$. To see that this
choice is measurable, we approximate it by continuous functions. For $\lambda>0$, let
\[
 z_\lambda(x)\coloneqq\argmin_{z\in\R^d}
       \left\{\bar f(z)+\frac{\|z-x\|^2}{2\lambda}\right\},
 \qquad h_\lambda(x)\coloneqq\frac{x-z_\lambda(x)}{\lambda}.
\]
The proximal minimizer is unique and depends continuously on $x$, and its optimality condition gives
$h_\lambda(x)\in\partial\bar f(z_\lambda(x))$. By monotonicity of the subdifferential,
$\|h_\lambda(x)\|^2\le\langle v(x),h_\lambda(x)\rangle$, so $\|h_\lambda(x)\|\le\|v(x)\|$. As
$\lambda\downarrow0$, $z_\lambda(x)\to x$, and the closed graph of the subdifferential places every limit
point of $h_\lambda(x)$ in $\partial\bar f(x)$. Uniqueness of the minimum-norm element gives
$h_\lambda(x)\to v(x)$. Thus $v$ is the pointwise limit of the continuous maps $h_{1/j}$ as the positive
integer $j$ tends to infinity, and its restriction to $\cX$ is the required bounded Borel selection.

\subsection{The contribution of each conditional variance}
\label{app:variance}

For every $t\in\{1,\ldots,n-1\}$, let $v_t\ge0$ be a finite deterministic bound on
$\E[\|W_t\|^2\mid\cF_t]$. The squared-distance inequality at time $t$ involves $X_{t+1}$. To include the last objective value
$f(X_n)$ in the same sum, take one auxiliary update from $X_n$ with $G_n\coloneqq v(X_n)$ and
$W_n\coloneqq0$. Let $M_n\coloneqq L^2$ and, for $1\le t<n$, let $M_t\coloneqq
L^2+v_t$. Since $X_1$ is deterministic and each noise has a finite second moment, the projected update
gives square-integrable iterates by induction. For $1\le k\le t\le n$, let $u\in\cX$ be square-integrable
and $\cF_k$-measurable. Given $\cF_t$, both $X_t-u$ and $G_t$ are known, so the cross terms involving
$W_t$ vanish in expectation. Expanding the squared distance after projection and using the subgradient
inequality therefore gives
\begin{equation}
 2\eta\,\E[f(X_t)-f(u)]
 \le \E[\|X_t-u\|^2]-\E[\|X_{t+1}-u\|^2]+\eta^2M_t .
 \label{eq:distance-progress}
\end{equation}
To recover the last objective value from the average controlled by this inequality, we drop one iterate
at a time. For $1\le k\le n$, let $a_k\coloneqq\E[f(X_k)]$ and let
$A_k\coloneqq(n-k+1)^{-1}\sum_{t=k}^n a_t$. The first suffix mean $A_1$ averages the entire run and the
last one is $A_n=a_n$. Taking $u=X_k$ in \eqref{eq:distance-progress} and summing over
this suffix gives $\sum_{t=k}^n(a_t-a_k)\le(\eta/2)\sum_{t=k}^n M_t$. Hence, for $1\le k<n$,
\[
 A_{k+1}-A_k
 =\frac{\sum_{t=k}^n(a_t-a_k)}{(n-k)(n-k+1)}
 \le\frac{\eta\sum_{t=k}^n M_t}{2(n-k)(n-k+1)}.
\]
For a deterministic $u\in\cX$, summing \eqref{eq:distance-progress} from $1$ to $n$ bounds $A_1-f(u)$ by
$\|X_1-u\|^2/(2\eta n)+(\eta/(2n))\sum_{t=1}^n M_t$. Now write $a_n=A_1+\sum_{k=1}^{n-1}(A_{k+1}-A_k)$. For each $1\le t\le n$, update $t$ appears in the
suffix starting at $k$ exactly when $k\le t$, so the coefficient of $\eta M_t/2$ is
\[
 \frac1n+\sum_{k=1}^{\min\{t,n-1\}}\frac1{(n-k)(n-k+1)}
 =\begin{cases}1/(n-t),&t<n,\\1,&t=n.\end{cases}
\]
Taking $u$ first as a nearest minimizer and then as $X_1$, and using $\Delta(X_1)\le LD_1$, yields
\begin{equation}
 \E[\Delta(X_n)]\le
 \min\left\{LD_1,\frac{D_1^2}{2\eta n}\right\}
 +\frac{\eta}{2}\left[L^2(1+H_{n-1})+\sum_{t=1}^{n-1}\frac{v_t}{n-t}\right].
 \label{gd:eq:heterogeneous-variance}
\end{equation}

For $2\le p\le\infty$, \eqref{eq:conditional-moment} permits $v_t=\sigma^2$ at every update. Thus the
last term in \eqref{gd:eq:heterogeneous-variance} is $(\eta/2)[L^2(1+H_{n-1})+\sigma^2 H_{n-1}]$. On the
normalized class we obtain
\begin{equation}
 \mathcal R_{n,p}^{(d)}(\eta,D)
 \le\min\left\{D,\frac{D^2}{2\eta n}\right\}+\eta(H_{n-1}+1/2).
 \label{gd:eq:variance-normalized}
\end{equation}
For $D>0$, choosing $\eta=D/\sqrt{n(2H_{n-1}+1)}$ gives the upper bound $D\sqrt{(2H_{n-1}+1)/n}$. In Appendix~\ref{sec:upper-proofs} we will start a scalar comparison at an iterate of small expected
error. Let $k_\star$ minimize $\E[f(X_k)]$ over $1\le k\le n$. Since its expected
error is at most $A_1-f_\star$, the same distance calculation with a nearest minimizer gives
\begin{equation}
 \E[\Delta(X_{k_\star})]\le\frac{D_1^2}{2\eta n}+\frac{\eta}{2}(L^2+\sigma^2).
 \label{eq:best-iterate}
\end{equation}
For $n=1$, Lipschitz continuity and the arithmetic-geometric mean inequality give $\Delta(X_1)\le
LD_1\le D_1^2/(2\eta)+\eta L^2/2$, so \eqref{eq:best-iterate} holds with $k_\star=1$ as well.

\subsection{A comparison for moments of order at most two}

Let $1\le p\le2$. We now want a noise vector with a finite second moment whose trajectory stays close to
the original one. Truncation can introduce a conditional bias, and subtracting that bias at every update
can move the iterate to a point with a different subgradient. We instead cancel the bias with a
correction made only on a small event. Until the first change in the noise, the two trajectories can
then use identical subgradients.

Let $\cG$ be a $\sigma$-field, let $W$ be an integrable random vector with $\E[W\mid\cG]=0$ and
$\E[\|W\|^p\mid\cG]\le\sigma^p$, and assume $\sigma>0$. For an event $E$, write $\1_E$ for its indicator.
Given a threshold $b\ge\sigma$, let $A\coloneqq W\1_{\{\|W\|\le b\}}$ and $\mu\coloneqq\E[A\mid\cG]$.
Since $\mu=-\E[W\1_{\{\|W\|>b\}}\mid\cG]$, we have $\|\mu\|\le\sigma^p b^{1-p}\le b$. Enlarge the
probability space by a uniform random variable $U$ on $[0,1]$ independent of the original variables. On
$\{\mu=0\}$ let $R\coloneqq0$, and on $\{\mu\ne0\}$ let
$R\coloneqq-b(\mu/\|\mu\|)\1_{\{U\le\|\mu\|/b\}}$. Given $\cG$, $A$ and $R$ are independent and
$\E[R\mid\cG]=-\mu$. The modified vector $\widehat W\coloneqq A+R$ is therefore centered given $\cG$ and
has norm at most $2b$. It satisfies
\begin{equation}
 \begin{aligned}
 \Pp[\widehat W\ne W\mid\cG]&\le2\sigma^p b^{-p},&
 \E[\|\widehat W-W\|\mid\cG]&\le2\sigma^p b^{1-p},\\
 \E[\|\widehat W\|^2\mid\cG]&\le2\sigma^p b^{2-p}.
 \end{aligned}
 \label{gd:eq:rare-change}
\end{equation}
For the first inequality, a change requires $\|W\|>b$ or a nonzero correction, whose conditional
probabilities are at most $\sigma^p b^{-p}$ and $\|\mu\|/b$. For the second, the expected discarded norm
is at most $\sigma^p b^{1-p}$ and the expected correction norm is $\|\mu\|$. For the third, conditional
independence gives the second moment as $\E[\|A\|^2\mid\cG]+b\|\mu\|-2\|\mu\|^2$, and each of the first
two terms is at most $\sigma^p b^{2-p}$.

For $1\le t<n$, choose a deterministic $b_t\ge\sigma$ and apply this construction to $W_t$ with a
fresh independent uniform $U_t$. Take all uniforms independent of the original process. For $1\le k\le
n$, enlarge $\cF_k$ to $\widehat\cF_k\coloneqq\cF_k\vee\sigma(U_1,\ldots,U_{k-1})$, with trivial added
$\sigma$-field at $k=1$. Independence gives the original conditional moments given $\widehat\cF_t$, and
\eqref{gd:eq:rare-change} applies to the resulting $\widehat W_t$. Start the comparison process at
$Y_1\coloneqq X_1$. At time $t$, if $Y_t=X_t$ use $\widehat G_t\coloneqq G_t$, and otherwise use
$\widehat G_t\coloneqq v(Y_t)$. Then let $Y_{t+1}\coloneqq\Pi_{\cX}(Y_t-\eta(\widehat G_t+\widehat
W_t))$. The subgradient choice is $\widehat\cF_t$-measurable, and the update keeps the trajectories
identical until a noise vector changes.

Nonexpansiveness of projection bounds the increase in $\|X_t-Y_t\|$ at update $t$ by $\eta\|W_t-\widehat
W_t\|+2\eta L\1_{\{X_t\ne Y_t\}}$. A modification at time $t$ can affect subgradient choices only in the
subsequent $n-t-1$ updates. Summing the distance inequalities, using $\1_{\{X_t\ne
Y_t\}}\le\sum_{s=1}^{t-1}\1_{\{W_s\ne\widehat W_s\}}$, and applying \eqref{gd:eq:rare-change} yield
\[
 \E[\|X_n-Y_n\|]\le\eta\sum_{t=1}^{n-1}
 \bigl[2\sigma^p b_t^{1-p}+4L(n-t-1)\sigma^p b_t^{-p}\bigr].
\]

We now choose the thresholds. An early modification can influence more subsequent subgradients, while
its variance carries a smaller weight in \eqref{gd:eq:heterogeneous-variance}. Both considerations favor
a larger threshold early in the run. Let $B\coloneqq\max\{L,\sigma\}$ and for $1\le t<n$ choose
$b_t\coloneqq B(n-t)$, so both the variance cost and the cost of changing the trajectory will sum powers
$(n-t)^{1-p}$. The conditional
second moment of $\widehat W_t$ is at most $2\sigma^p B^{2-p}(n-t)^{2-p}$. In
\eqref{gd:eq:heterogeneous-variance} the resulting noise contribution for $Y_n$ is at most $\eta\sigma^p
B^{2-p}S_{n-1,p}$. By Lipschitz continuity, passing from $Y_n$ to $X_n$ adds at most $L\E[\|X_n-Y_n\|]$.
Using $n-t-1\le n-t$ in the preceding display, the combined coefficient of $\eta\sigma^p B^{2-p}S_{n-1,p}$
is at most $1+2L/B+4(L/B)^2\le7$. Consequently,
\begin{equation}
 \E[\Delta(X_n)]\le\min\left\{LD_1,\frac{D_1^2}{2\eta n}\right\}
 +\frac{\eta L^2}{2}(1+H_{n-1})+7\eta\sigma^p B^{2-p}S_{n-1,p}.
 \label{gd:eq:general-moment-upper}
\end{equation}
When $\sigma=0$, \eqref{gd:eq:heterogeneous-variance} gives the same bound with the last term omitted.
When $L=0$, the error vanishes.

For $1\le p\le2$, the normalized class has $L,\sigma,B\le1$ and $S_{n-1,p}\ge\max\{1,H_{n-1}\}$. Substituting
these bounds and \eqref{eq:sum-comparison} into \eqref{gd:eq:general-moment-upper} gives
\begin{equation}
 \mathcal R_{n,p}^{(d)}(\eta,D)
 \le\min\left\{D,\frac{D^2}{2\eta n}\right\}+8\eta\Phi_{n,p}.
 \label{gd:eq:heavy-fixed-step}
\end{equation}

\section{Upper bounds on the line}
\label{sec:upper}

Throughout this appendix $d=1$. To control how much the objective can increase after a given iterate, we
fix a sublevel interval and measure the distance outside its two ends. Beyond either endpoint, convexity makes the noiseless update point back toward the interval. We first translate this fact into an inequality for the distance, then bound the
expected objective increase during the excursion. Finally, we choose the starting level according to
whether the noise has a finite second moment.

\subsection{Distances outside a sublevel interval}
\label{sec:geometry}

For a real number $v$, write $v_+\coloneqq\max\{v,0\}$. Let $b\in(0,\infty]$, let $I$ be $[0,b]$ if
$b<\infty$ and $[0,\infty)$ otherwise, and let $m:I\to[0,L]$ be nondecreasing. Think of $z$ as the distance beyond one endpoint and $m(z)$ as the slope there. The noiseless step
leaves at most $(z-\eta m(z))_+$ of that distance, but this quantity can decrease at a jump of the
slope. To compare trajectories started at different distances, we take its nondecreasing upper envelope.
Following Beretta et al.~\cite[Definition~4 and Proposition~5]{BCCP2026}, for $z\in I$ define
\[
 \mathcal T(z)\coloneqq\sup_{0\le s\le z}(s-\eta m(s))_+,\qquad
 J(z)\coloneqq\int_0^z(s-\mathcal T(s))\,ds .
\]
As the set over which we take the supremum grows with $z$, $\mathcal T$ is nondecreasing. Every term is
at most $z$ and the choice $s=z$ gives $\mathcal T(z)\ge(z-\eta m(z))_+$, hence $0\le z-\mathcal
T(z)\le\eta m(z)$. For $x,y\in I$ with $x\le y$, terms with $s\le x$ in the supremum for $\mathcal T(y)$
are at most $\mathcal T(x)$, while for $x<s\le y$, monotonicity gives $(s-\eta m(s))_+\le(x-\eta
m(x))_++s-x$. Thus $0\le\mathcal T(y)-\mathcal T(x)\le y-x$, which makes both $\mathcal T$ and $z\mapsto
z-\mathcal T(z)$ nondecreasing and $1$-Lipschitz. For finite $b$ and $z>b$, extend $\mathcal T$ by
$\mathcal T(z)\coloneqq z-b+\mathcal T(b)$ and use the same integral to extend $J$. The stated
monotonicity and Lipschitz properties persist, $\mathcal T$ is unbounded, and $J$ is nondecreasing and
$\eta L$-Lipschitz. When $b=\infty$, unboundedness follows from $\mathcal T(z)\ge z-\eta L$.

We now relate the potential $J$ to the increase of the objective along this branch. For $z\in I$, the
increase obtained by integrating the slope is $h(z)\coloneqq\int_0^z m(s)\,ds$, and we claim that
\begin{equation}
 h(z)\le J(z)/\eta+(\eta L^2/2)\1_{\{z>0\}} .
 \label{eq:h-by-envelope}
\end{equation}
To compare the two integrals, fix a finite $z\in I$ with $z>0$ and compare how far along $[0,z]$ each integrand first
reaches the same slope level. For $u\in[0,L]$, let $a(u)$ and $c(u)$ be the infima in $[0,z]$
where $m(s)\ge u$ and $s-\mathcal T(s)\ge\eta u$, respectively, assigning the value $z$ to an empty set.
Since $s-\mathcal T(s)\le\eta m(s)$, we have $a(u)\le c(u)$. If $a(u)+\eta u\le z$, then for $s\le a(u)$
we have $(s-\eta m(s))_+\le a(u)$, and for $a(u)<s\le a(u)+\eta u$ monotonicity gives $m(s)\ge u$. Hence
$\mathcal T(a(u)+\eta u)\le a(u)$. In this case $c(u)\le a(u)+\eta u$, and the inequality also holds if
$a(u)+\eta u>z$. Integrating the lengths of the two superlevel intervals gives
$h(z)-J(z)/\eta=\int_0^L(c(u)-a(u))\,du\le\eta L^2/2$. At zero both sides of \eqref{eq:h-by-envelope}
vanish.

For a nonempty set $S\subseteq\R$ and $x\in\R$, write $\dist(x,S)\coloneqq\inf_{y\in S}|x-y|$. Fix $a\ge
f_\star$. Let $S_a\coloneqq\{x\in\cX:f(x)\le a\}$ and $E_a\coloneqq\{x\in\cX:\dist(x,S_a)\le\eta L\}$.
Both are nonempty closed intervals. The enlargement allows for a whole noiseless step from $S_a$, whose
length is at most $\eta L$. Within this enlarged interval the error above $a$ is already at most $\eta
L^2$, so it suffices to control any distance beyond its endpoints. We use a sign $\varepsilon\in\{-1,1\}$ to treat the two ends of $E_a$
together. Let $\mathcal B_a$ be the set of signs for which $\sup_{x\in E_a}\varepsilon
x<\sup_{x\in\cX}\varepsilon x$. For each $\varepsilon\in\mathcal B_a$, let
$b_\varepsilon\coloneqq\sup_{x\in E_a}\varepsilon x$ and $d_\varepsilon\coloneqq\sup_{x\in\cX}\varepsilon
x-b_\varepsilon$. For $x\in\cX$ and finite $0\le z\le d_\varepsilon$, define the distance and objective
increase on that side by
\[
 z_\varepsilon(x)\coloneqq(\varepsilon x-b_\varepsilon)_+,\qquad
 h_\varepsilon(z)\coloneqq f(\varepsilon(b_\varepsilon+z))-f(\varepsilon b_\varepsilon).
\]
The branch function $h_\varepsilon$ is convex, nondecreasing, $L$-Lipschitz and zero at the origin. Its
left derivative at positive $z$, extended by zero at $z=0$, is a nondecreasing function $m_\varepsilon$
with integral $h_\varepsilon$. Denote the resulting map and potential by $\mathcal T_\varepsilon$ and
$J_\varepsilon$. This enlargement of the sublevel interval follows Beretta et al.~\cite[Lemma~10 and
Appendix~A.3]{BCCP2026}.

For $1\le j<n$ and $\varepsilon\in\mathcal B_a$, the update satisfies
\begin{equation}
 z_\varepsilon(X_{j+1})
 \le\bigl(\mathcal T_\varepsilon(z_\varepsilon(X_j))-\varepsilon\eta W_j\bigr)_+ .
 \label{eq:branch-recursions}
\end{equation}
To see this, consider the right side of the interval, for which $\varepsilon=1$ and $b_1=\sup S_a+\eta
L$. Let $x\in\cX$ and let $g\in\partial_{\cX}f(x)$ have $|g|\le L$. If $x\le b_1-\eta L$, then $x-\eta
g\le b_1$. If $b_1-\eta L<x\le b_1$, convexity gives $g\ge0$ and again $x-\eta g\le b_1$. If $x>b_1$, the
subgradient is at least the left derivative at $x$, giving $x-\eta g-b_1\le\mathcal T_1(z_1(x))$. Since
projection onto an interval containing $b_1$ cannot increase the positive excess above $b_1$, adding the
noise and projecting gives \eqref{eq:branch-recursions}. Reflecting the line proves the left-side
inequality. For $x\in E_a$, Lipschitz continuity gives $(f(x)-a)_+\le\eta L^2$. At every $x\in\cX$ at
most one branch has positive distance, so adding its excess from \eqref{eq:h-by-envelope} gives
\begin{equation}
 (f(x)-a)_+\le\tfrac32\eta L^2+\frac1\eta
      \sum_{\varepsilon\in\mathcal B_a}J_\varepsilon(z_\varepsilon(x)).
 \label{eq:branch-potential}
\end{equation}

\subsection{A uniform bound for a reflected recursion}
\label{sec:horizon}

We now bound the potential for either branch. The increments may depend on the history. It will suffice
to dominate their expected excess above each threshold by that of one reference variable. The next lemma
turns this one-step comparison into a bound valid at every future time.

\begin{lemma}\label{thm:horizon-uniform}
Let $\mathcal T:[0,\infty)\to[0,\infty)$ be unbounded, nondecreasing, and $1$-Lipschitz, and assume that
for every $z\ge0$, $\mathcal T(z)\le z$. Let $(\cF_t)_{t\ge1}$ be a filtration and let $(Z_t)_{t\ge1}$ be
a nonnegative adapted process starting at zero. For each $t\ge1$, let $R_t$ be a real
$\cF_{t+1}$-measurable increment satisfying $Z_{t+1}\le(\mathcal T(Z_t)+R_t)_+$. Let $V$ be a
nondegenerate square-integrable reference variable with an even density $v$ that is completely monotone
on $(0,\infty)$, meaning that for every integer $j\ge0$ and every $x>0$ its derivative satisfies
$(-1)^jv^{(j)}(x)\ge0$. Suppose that on one event of probability one, for every $t\ge1$ and $s\in\R$,
$\E[(R_t-s)_+\mid\cF_t]\le\E[(V-s)_+]$. Then
\begin{equation}
 \sup_{t\ge1}\E\left[\int_0^{Z_t}(z-\mathcal T(z))\,dz\right]
 \le\tfrac12\E[V^2].
 \label{eq:horizon-bound}
\end{equation}
More generally, if the recursion starts at a deterministic time $k$ and $Z_k=0$ on an event $A\in\cF_k$,
then for every $t\ge k$ the same bound holds with $\1_A$ inside the expectation on the left and the
factor $\Pp[A]$ on the right.
\end{lemma}

\begin{proof}
The assumption controls positive parts of increments, so we first express the potential in terms of
positive parts of the state. We will then bound those positive parts recursively and sum the bounds
along a reference walk. For $s\ge0$, let $I(s)\coloneqq\inf\{y\ge0:\mathcal T(y)\ge s\}$, and for $y\ge0$ let $B(y)\coloneqq
y-\mathcal T(y)$. By continuity, monotonicity and unboundedness, $\mathcal T(I(s))=s$ and $I(s)\ge s$.
Since $B$ is nondecreasing, $(\mathcal T(y)-s)_+\le(y-I(s))_+$. The set $\{s\ge0:I(s)\le y\}$ equals
$[0,\mathcal T(y)]$, so the image of Lebesgue measure under $I$ is the Stieltjes measure $d\mathcal T$.
Consequently, for every nonnegative nonincreasing integrable $\varphi$ on $[0,\infty)$ and every $z\ge0$,
\[
 \int_0^\infty\bigl(\varphi(s)-\varphi(I(s))\bigr)\,ds
   =\int_{[0,\infty)}\varphi(y)\,dB(y),\qquad
 \int_0^z B(s)\,ds=\int_{[0,\infty)}(z-y)_+\,dB(y).
\]
The second identity follows by integrating $B(s)=\int_{[0,s]}dB(y)$.

The second identity reduces our task to bounding $\E[(Z_t-x)_+]$ over all thresholds $x\ge0$. To obtain
a deterministic upper bound, we use an independent draw from the law of $V$ at each recursion step. For
$x\ge0$, let $g(x)\coloneqq\E[(V-x)_+]$, let $c_1(x)\coloneqq0$, and for $t\ge1$ let
$c_{t+1}(x)\coloneqq\E[c_t(I((x-V)_+))]+g(x)$. We claim that $\E[(Z_t-x)_+]\le c_t(x)$. The claim starts
at $t=1$. The comparison at threshold zero makes $(R_t)_+$ integrable, hence each $Z_t$ is integrable.
Since positive parts are $1$-Lipschitz in the threshold, approximation by finite-valued rational
thresholds extends the conditional comparison to the integrable random threshold $x-\mathcal T(Z_t)$. For
$w\in\R$ and $y\ge0$ we have
\[
 (\mathcal T(y)+w-x)_+
 =(\mathcal T(y)-(x-w)_+)_++(w-x)_+
 \le(y-I((x-w)_+))_++(w-x)_+ .
\]
Since $x\ge0$, the recursion bounds $(Z_{t+1}-x)_+$ by $(\mathcal T(Z_t)+R_t-x)_+$. At the threshold
$x-\mathcal T(Z_t)$, the conditional comparison bounds its expectation by that for an independent copy
of $V$ in place of $R_t$. For each value $w$ of that copy, the induction hypothesis applies at the
deterministic threshold $I((x-w)_+)$. Integrating the two terms of the displayed inequality over $w$
gives $c_{t+1}(x)$, as required. For a start at $k$ on $A\in\cF_k$, we can set the state and all subsequent increments to zero off $A$.
The comparison there follows from $(-s)_+\le\E[(V-s)_+]$ by Jensen's inequality. On $A$ the same
induction, with each term multiplied by $\1_A$, gives
$\E[\1_A(Z_t-x)_+]\le\Pp[A]c_{t-k+1}(x)$.

The recursion and the monotonicity of $I$ give by induction that each $c_t$ is finite, nonnegative, and
nonincreasing in $x$, and that $c_{t+1}\ge c_t$. To check integrability, note that $\int_0^\infty
g(x)\,dx=\E[(V_+)^2]/2$, and for any finite nonnegative nonincreasing integrable $h$, splitting the
integral according to the sign of $V$ gives
$\int_0^\infty\E[h((x-V)_+)]\,dx\le\int_0^\infty h(x)\,dx+\E[V_+]h(0)$. Since $I(s)\ge s$, this proves
integrability by induction.

We now need a bound uniform in $t$. For $\psi_t\coloneqq c_t-c_t\circ I$, the first Stieltjes identity
gives $\int_{[0,\infty)}c_t(y)\,dB(y)=\int_0^\infty\psi_t(x)\,dx$. The recursion for $c_{t+1}$ evaluates
$c_t\circ I$ at $(x-V)_+$ and adds $g(x)$. Since an excursion of the walk with this reflected update starts and ends at zero,
summing the recursion over one excursion will bound the expected sum of $\psi_t$ over the heights it
visits. To turn this estimate into a bound on the ordinary integral of $\psi_t$, we will use the
complete monotonicity of the reference density to show that the expected number of visits per unit of
height has a positive lower bound throughout $(0,\infty)$.

Let $(V_j)_{j\ge1}$ be independent copies of $V$. Let $S_0\coloneqq\zeta_0\coloneqq0$ and, for $j\ge1$,
let $S_j\coloneqq\sum_{i=1}^j V_i$ and $\zeta_j\coloneqq(\zeta_{j-1}-V_j)_+$. Write $\E_0$ for
expectation on this space. Let $\tau^-$ and $\tau^+$ be the first times at which $S$ is strictly negative and strictly positive, respectively, and let $H^-\coloneqq-S_{\tau^-}$ and $H^+\coloneqq S_{\tau^+}$ be the
corresponding distances beyond zero. The centered nondegenerate walk oscillates, so both times are
finite almost surely~\cite[Chapter~XII]{Feller1971}. Since the increments have a density, almost surely
$S$ never returns exactly to zero at a positive time. For $0\le j<\tau^+$ we therefore have $\zeta_j=-S_j$, and at $\tau^+$ the reflected walk first returns
to zero.

For a nonnegative Borel function $h$ on $[0,\infty)$, let $Q(h)\coloneqq\E_0[\sum_{j=0}^{\tau^+-1}
h(\zeta_j)]$. Conditional on the first $j$ increments, $g(\zeta_j)$ is the expected value of
$(V_{j+1}-\zeta_j)_+$. During the excursion this positive part vanishes at every step except the return,
where it equals $H^+$. Consequently $Q(g)=\E[H^+]$. Doney's ladder-height
identity~\cite[p.~248]{Doney1980} gives finite means and $\E[V^2]=2\E[H^-]\E[H^+]$, with absolute
continuity making the strict and weak crossings equivalent.

Since $\zeta_0=\zeta_{\tau^+}=0$, summing a function over $\zeta_1,\ldots,\zeta_{\tau^+}$ gives the same
value as summing it over $\zeta_0,\ldots,\zeta_{\tau^+-1}$. Using this shift in the recursion for
$c_{t+1}$ gives $Q(c_{t+1})=Q(c_t\circ I)+\E[H^+]$. Since $I(x)\ge x$ and $c_t$ is nonincreasing,
induction from $Q(c_1)=0$ gives finiteness. Subtracting the two finite expectations and using
$c_{t+1}\ge c_t$ now yields
\[
 Q(\psi_t)=\E[H^+]-\bigl(Q(c_{t+1})-Q(c_t)\bigr)\le\E[H^+] .
\]

It remains to relate this sum over visited heights to the integral of $\psi_t$. For each integer $j\ge1$, the excursion lasts past $j$ precisely on $\{S_1<0,\ldots,S_j<0\}$, up to a
null event. On this event $\zeta_j=-S_j$. If we
reverse the first $j$ increments, the earlier partial sums of the reversed path are $S_j-S_i$ for $1\le
i\le j$. The original partial sums are all negative precisely when the reversed path ends at a new
strict minimum. Reversal preserves the law, so summing over $j$ counts the successive strict descending
ladder heights of $S$, that is, the distances below zero at successive new minima. By independence of
the increments after each such time, the increases in those distances are independent copies of $H^-$.

We now use complete monotonicity to bound the density of these visits from below. To find the density of $H^-$ at $x>0$, sum over the possible crossing times $j$. Given the path through
time $j-1$ with no earlier descent, the crossing step must take $S_{j-1}$ to $-x$, so symmetry of $v$
gives
\[
 k(x)\coloneqq\sum_{j=1}^\infty\E\bigl[\1_{\{\tau^->j-1\}}v(x+S_{j-1})\bigr].
\]
On each event in this sum $S_{j-1}\ge0$. By Bernstein's representation
theorem~\cite[Theorem~1.4]{Schilling2012}, $v$ on $(0,\infty)$ is a positive mixture of the functions
$x\mapsto e^{-rx}$ for $r\ge0$. Translating by $S_{j-1}$ only multiplies each integrand by
$e^{-rS_{j-1}}$, so Tonelli's theorem gives the same kind of representation for $k$. Since $k$ integrates to one, it is finite on $(0,\infty)$ and completely monotone there. Moreover, if
$\nu$ is its representing measure, then $\nu(\{0\})=0$ and $\int_{(0,\infty)}r^{-1}\,\nu(dr)=1$. Thus
$k$ is a mixture of the exponential densities $x\mapsto re^{-rx}$ with probability measure
$r^{-1}\nu(dr)$.

For $j\ge1$ let $k^{*j}$ denote the $j$-fold convolution of $k$. The successive descending ladder
heights therefore have total density $u(x)\coloneqq\sum_{j=1}^\infty k^{*j}(x)$ for $x>0$. The renewal
representation for mixtures of exponentials~\cite[Theorem~2.1]{EnriquezNoiry2020} gives a nonincreasing
version of $u$, and for every $h>0$, Blackwell's renewal theorem~\cite{Blackwell1953} gives
$\int_x^{x+h}u(y)\,dy\to h/\E[H^-]$ as $x\to\infty$. Thus $u(x)\ge1/\E[H^-]$ throughout $(0,\infty)$.

The excursion visits zero at its start. Counting that visit and then the descending ladder heights gives
$Q(\varphi)=\varphi(0)+\int_0^\infty\varphi(x)u(x)\,dx$ for every nonnegative Borel function $\varphi$.
Since $I(x)\ge x$ and $c_t$ is nonincreasing, $\psi_t\ge0$ and the density bound gives $\int_0^\infty\psi_t(x)\,dx\le\E[H^-]Q(\psi_t)$. Returning
to the two Stieltjes identities and the positive-part comparison, we obtain
\[
 \E\left[\int_0^{Z_t}B(s)\,ds\right]
 \le\int_{[0,\infty)}c_t(y)\,dB(y)
 =\int_0^\infty\psi_t(x)\,dx
 \le\E[H^-]\E[H^+]=\tfrac12\E[V^2].
\]
For a start at time $k$ on $A$, use $\E[\1_A(Z_t-y)_+]\le\Pp[A]c_{t-k+1}(y)$ in the first inequality.
The same argument multiplies the right-hand side by $\Pp[A]$, as claimed.
\end{proof}

\subsection{Reference distributions for scalar noise}
\label{sec:reference}

To apply Lemma~\ref{thm:horizon-uniform}, we need a reference law with a tractable second moment. For
heavy-tailed noise, we first truncate at a level we can later choose. Let $\cG$ be a $\sigma$-field and
let $Q$ be an integrable real variable with $\E[Q\mid\cG]=0$. For $K>0$, we construct a centered
modification $\widehat Q$ on the same space with $|\widehat Q|\le\min\{|Q|,K\}$. Clipping may leave the
positive and negative parts with different conditional means, so we scale down the part with the larger
mean. Let $Q_K\coloneqq\max\{-K,\min\{Q,K\}\}$ and let $P\coloneqq(Q_K)_+$,
$N\coloneqq(-Q_K)_+$, $\mu_+\coloneqq\E[P\mid\cG]$, and $\mu_-\coloneqq\E[N\mid\cG]$. If $\mu_+\ge\mu_-$,
let $\widehat Q\coloneqq(\mu_-/\mu_+)P-N$, interpreting $0/0$ as one, and if $\mu_+<\mu_-$, let $\widehat
Q\coloneqq P-(\mu_+/\mu_-)N$. Each multiplier lies in $[0,1]$, and the scaled positive and negative parts
now have equal conditional means. Since $\E[|Q_K-\widehat
Q|\mid\cG]=|\mu_+-\mu_-|\le\E[(|Q|-K)_+\mid\cG]$, the triangle inequality gives
\[
 \E[|Q-\widehat Q|\mid\cG]\le2\E[(|Q|-K)_+\mid\cG].
\]
In particular, if $1\le p<\infty$, $\sigma>0$, $x\ge1$, $\E[|Q|^p\mid\cG]\le\sigma^p$, and $K=x\sigma$,
then this bound is at most $2\sigma x^{1-p}$.

Let $1\le p<\infty$, $\sigma>0$, and $x\ge1$, and assume $\E[|Q|^p\mid\cG]\le\sigma^p$ and $|Q|\le
x\sigma$. We now choose the reference law for the bounded variable. Its magnitude must cover the tail allowed by
the moment bound, and the reflected-recursion lemma requires an even density that is completely monotone
on the positive half-line. Take independent variables $\varepsilon,E,B$, where $\varepsilon$ is a
uniform sign and $E$ is exponential with mean one. To implement the tail comparison, for $y\ge0$ choose
$B$ with
\[
 \Pp[B>y]\coloneqq
 \begin{cases}
  1,&0\le y<2\sigma,\\
  (2\sigma/y)^p,&2\sigma\le y<2x\sigma,\\
  0,&y\ge2x\sigma.
 \end{cases}
\]
Let $V\coloneqq\varepsilon BE$. For $u\ge0$, conditional Markov's inequality and the cutoff give
$\Pp[2|Q|>u\mid\cG]\le\Pp[B>u]$. For $s\ge0$, tail integration and $(Q-s)_+\le(2|Q|-s)_+/2$ give
$\E[(Q-s)_+\mid\cG]\le\E[(\varepsilon B-s)_+]$. For a centered $Q$, the identity $\E[(Q-s)_+\mid\cG]=-s+\E[(-Q+s)_+\mid\cG]$ reduces a negative
threshold to the positive threshold $-s$ for $-Q$. Symmetry of the reference law then gives the same
comparison. The exponential factor supplies the required density. Since $\E[V\mid\varepsilon B]=\varepsilon B$,
conditional Jensen's inequality gives $\E[(\varepsilon B-s)_+]\le\E[(V-s)_+]$ at every threshold $s$. The positive-part expectations are
$1$-Lipschitz in $s$, so applying the comparison first on the rationals gives versions for which it holds
simultaneously at all real thresholds. The density of $V$ is even, and on the positive half-line it
equals $t\mapsto\E[B^{-1}\exp(-t/B)]/2$, a mixture of decreasing exponentials. Since $\E[E^2]=2$, tail
integration gives
\begin{equation}
 \tfrac12\E[V^2]=\E[B^2]
 =2\int_0^{2x\sigma}y\Pp[B>y]\,dy
 =4\sigma^2+8\sigma^2\int_1^x u^{1-p}\,du.
 \label{eq:reference-second-moment}
\end{equation}
If $|Q|\le\sigma$, then for every $s\in\R$ the chord of $q\mapsto(q-s)_+$ on $[-\sigma,\sigma]$ and
conditional centering give the comparison with $\varepsilon\sigma$. Jensen's inequality then gives the
reference variable $\varepsilon\sigma E$, with $\E[(\varepsilon\sigma E)^2]/2=\sigma^2$.

\subsection{From a sublevel interval to the last iterate}
\label{sec:upper-proofs}

We first start from a deterministic level $a$ at a time when we know that the iterate lies below it.
This lets us start both exterior distances at zero, as required by Lemma~\ref{thm:horizon-uniform}. Fix
$a\ge f_\star$, $1\le k\le n$, and $A\in\cF_k$ such that $f(X_k)\le a$ on $A$. Let
$r\coloneqq|\mathcal B_a|$ be the number of exterior branches. For finite $p$ and $x\ge1$, write
$I_p(x)\coloneqq\int_1^x u^{1-p}\,du$. We will use the bound
\begin{equation}
 \E[\1_A(f(X_n)-a)_+]\le
 \eta\Pp[A]\left[\tfrac32 L^2+r\sigma^2(4+8I_p(x))
                       +2L(n-k)\sigma x^{1-p}\right].
 \label{eq:event-terminal}
\end{equation}
If the noises from time $k$ onward already satisfy $|W_j|\le x\sigma$, the final term can be dropped. For
$p=\infty$, the bracket becomes $\tfrac32 L^2+r\sigma^2$.

On $A$ each branch coordinate starts at zero at time $k$. For finite $p$, $\sigma>0$ and $k\le j<n$,
construct $\widehat W_j$ from $W_j$ with cutoff $x\sigma$ as above, using the same noise modifications
for both branches. For each $\varepsilon\in\mathcal B_a$, run the recursion in
\eqref{eq:branch-recursions} with equality, starting at zero, first with the noises $W_j$, then with
$\widehat W_j$. Write $Y_{\varepsilon,j}$ and $\widehat Y_{\varepsilon,j}$ for the resulting coordinates.
Monotonicity gives $z_\varepsilon(X_j)\le Y_{\varepsilon,j}$ on $A$, and the $1$-Lipschitz update gives
$|Y_{\varepsilon,n}-\widehat Y_{\varepsilon,n}|\le\eta\sum_{j=k}^{n-1}|W_j-\widehat W_j|$. By the localized
bound in Lemma~\ref{thm:horizon-uniform} and \eqref{eq:reference-second-moment} at scale $\eta\sigma$,
\[
 \E[\1_A J_\varepsilon(\widehat Y_{\varepsilon,n})]
 \le\Pp[A]\eta^2\sigma^2(4+8I_p(x)).
\]
A finite recursion may be extended by zero increments, which satisfy the comparison by centering of the
reference variable. At $X_n$ at most one actual branch coordinate is positive, and $J_\varepsilon$ is
$\eta L$-Lipschitz. Thus replacing the modified coordinates by the actual coordinates in the sum of
potentials costs at most $\eta^2 L\sum_{j=k}^{n-1}\E[\1_A|W_j-\widehat W_j|]$, even when both modified
coordinates are positive. Each summand is at most $2\Pp[A]\sigma x^{1-p}$ by conditional expectation.
Dividing by $\eta$ in \eqref{eq:branch-potential} proves \eqref{eq:event-terminal}. For noises already
within the cutoff the two processes coincide. For $p=\infty$, applying the comparison for bounded increments from
Appendix~\ref{sec:reference} at scale $\eta\sigma$ gives the stated bound, and for $\sigma=0$ each
branch starting at zero remains there.

We next want to choose the level after observing $X_k$. To use the deterministic-level estimate, we
round $\Delta(X_k)$ upward to a multiple of $\delta>0$ and apply the estimate separately on each event
of the resulting partition. The first event is $\{\Delta(X_k)=0\}$, and for integers $j\ge1$ the others are
$\{(j-1)\delta<\Delta(X_k)\le j\delta\}$. On each event use the corresponding
upper endpoint $f_\star+j\delta$ as the deterministic level, with $j=0$ for the first event. The positive
part changes by at most $\delta$. Summing the localized inequalities and letting $\delta\downarrow0$
gives the same bounds on $\E[(f(X_n)-f(X_k))_+]$ with probability factor one. We may use two branches in
general, and one if every sublevel interval has at most one exterior branch.

We now handle the initial distance. With a finite variance, \eqref{eq:best-iterate} supplies a low-error
iterate and the random-level bound controls what happens afterward. For $1\le p<2$ that variance bound
need not be available. We start each branch at its actual initial distance. After truncation the comparison recursion is itself
a subgradient method for its potential, which lets us bound that potential from an arbitrary starting
point. Let $N\ge2$
be an integer, let $b\ge0$, and let $\mathcal T$ be a nondecreasing $1$-Lipschitz map on $[0,\infty)$
such that, for every $z\ge0$, $0\le z-\mathcal T(z)\le b$ and $\mathcal T(z)\ge0$. Let $y_0\ge0$ be
deterministic, let $Z_1\coloneqq y_0$, let $(\cF_t)_{t=1}^N$ be a filtration, and fix $x\ge1$. For $1\le
t<N$, let $R_t$ be $\cF_{t+1}$-measurable, with $|R_t|\le x\sigma$, $\E[R_t\mid\cF_t]=0$ and
$\E[|R_t|^p\mid\cF_t]\le\sigma^p$, and let $Z_{t+1}\coloneqq(\mathcal T(Z_t)+R_t)_+$. The function
$J(z)\coloneqq\int_0^z(s-\mathcal T(s))\,ds$ is convex, nondecreasing and $b$-Lipschitz, and
$J'(z)=z-\mathcal T(z)$. Thus this recursion is projected Stochastic subGradient Method for $J$ on
$[0,\infty)$, with unit stepsize and noise $-R_t$. Since $|R_t|\le x\sigma$, its conditional noise
variance is at most $(x\sigma)^{2-p}\E[|R_t|^p\mid\cF_t]\le\sigma^2x^{2-p}$. Applying
\eqref{eq:distance-progress} with comparator zero and averaging gives an index $k$ with $\E[J(Z_k)]\le
y_0^2/(2N)+(b^2+\sigma^2x^{2-p})/2$. Since $J$ is nondecreasing, every sublevel interval has at most one
exterior branch. Adding the corresponding bounded-noise estimate for $\E[(J(Z_N)-J(Z_k))_+]$ gives
\begin{equation}
 \E[J(Z_N)]\le\frac{y_0^2}{2N}+2b^2+
       \sigma^2\bigl[4+\tfrac12x^{2-p}+8I_p(x)\bigr].
 \label{eq:arbitrary-initial}
\end{equation}

Fix a cutoff factor $x\ge1$. For $1\le p<2$, take $a=f_\star$. On each present branch
$\varepsilon\in\mathcal B_{f_\star}$, start the original-noise and modified-noise recursions at
$z_\varepsilon(X_1)$. Only the branch on the side containing $X_1$ can have positive initial distance. Since $E_{f_\star}$
contains every minimizer, the sum of the squared starting coordinates is $\dist(X_1,E_{f_\star})^2\le
D_1^2$.
Apply \eqref{eq:arbitrary-initial} to the modified recursions with $N=n$, $b=\eta L$, and scale
$\eta\sigma$. Passing back to $X_n$ again costs at most $2\eta^2L\sigma(n-1) x^{1-p}$ in the sum of
potentials. For $2\le p<\infty$, use the best index $k_\star$ from \eqref{eq:best-iterate} and bound
$\E[(f(X_n)-f(X_{k_\star}))_+]$ by the two-branch form of \eqref{eq:event-terminal}. For every $x\ge1$,
using $\Delta(X_n)\le\Delta(X_{k_\star})+(f(X_n)-f(X_{k_\star}))_+$ in the second case gives the two
bounds
\[
 \E[\Delta(X_n)]\le\frac{D_1^2}{2\eta n}+2\eta L\sigma(n-1) x^{1-p}
 +\eta
 \begin{cases}
  \tfrac{11}{2}L^2+\sigma^2(8+x^{2-p}+16I_p(x)),&1\le p<2,\\
  2L^2+\sigma^2(\tfrac{17}{2}+16I_p(x)),&2\le p<\infty .
 \end{cases}
\]
The clipping error contains $(n-1)x^{1-p}$, while keeping more of the noise increases $I_p(x)$. Taking $x=n-1$
gives $(n-1)x^{1-p}=(n-1)^{2-p}\le\Phi_{n,p}$. Using this and $I_p(n-1)\le\Phi_{n,p}$, and using the bounded reference law
directly at $p=\infty$, we obtain
\begin{equation}
 \E[\Delta(X_n)]\le\frac{D_1^2}{2\eta n}+\eta
 \begin{cases}
 \tfrac{11}{2}L^2+(25\sigma^2+2L\sigma)\Phi_{n,p},&1\le p<2,\\
 2L^2+(\tfrac{49}{2}\sigma^2+2L\sigma)\Phi_{n,p},&2\le p<\infty,\\
 2L^2+\tfrac52\sigma^2,&p=\infty.
 \end{cases}
 \label{eq:upper}
\end{equation}

Finally, starting \eqref{eq:event-terminal} at the deterministic level $a=f(X_1)$, using $\Delta(X_1)\le
LD_1$ and taking $x=n-1$ gives the initial term $LD_1$. For finite $p$ the accompanying term is at most
$\eta[\tfrac32L^2+(24\sigma^2+2L\sigma)\Phi_{n,p}]$, and for $p=\infty$ it is at most
$\eta(\tfrac32L^2+2\sigma^2)$. On the normalized class, combining this with \eqref{eq:upper} yields
\begin{equation}
 \mathcal R_{n,p}^{(1)}(\eta,D)
 \le\min\left\{D,\frac{D^2}{2\eta n}\right\}+\tfrac{65}{2}\eta\Phi_{n,p}.
 \label{eq:normalized-upper}
\end{equation}

\section{The upper bound above the second moment}
\label{app:geometric}

Fix $d\ge2$, $2\le p\le\infty$, $D\ge0$, and an instance in $\cC_p^{(d)}(D)$. In this appendix
$C$ denotes a positive universal constant that may change from one estimate to the next. We first take $\eta=1$ and subtract
$f_\star$ from the objective. Write $\xi_k\coloneqq-W_k$, so that the update adds $\xi_k$ after
subtracting $G_k$. Let $\varepsilon\coloneqq\min\{p-2,1\}$, with $\varepsilon=1$ at $p=\infty$.
For a random variable $U$, write $\E_s[U]\coloneqq\E[U\mid\cF_s]$ and, for an event $B$, write
$\Pp_s[B]\coloneqq\Pp[B\mid\cF_s]$. For $q\ge1$ we denote the conditional and unconditional moment
norms of a random vector $V$ by $\|V\|_{q,s}\coloneqq(\E_s[\|V\|^q])^{1/q}$ and
$\|V\|_q\coloneqq(\E[\|V\|^q])^{1/q}$, respectively. By the noise assumption and conditional Jensen's
inequality, $\E_k[\xi_k]=0$ and $\E_k[\|\xi_k\|^{2+\varepsilon}]\le1$ at each update.

The variance estimate \eqref{gd:eq:variance-normalized} can charge a constant amount of error to
each of order $\log n$ time scales. To see how many of these scales can contribute in dimension $d$,
we will compare displacements over time intervals ending at the same iterate. We first control the
motion within one interval. The point of this first step is to distinguish the contribution of large
noise increments, whose cost depends on $p$, from displacements whose fourth moments we can bound
uniformly. We then use convexity to relate displacements on different intervals and use the dimension
to bound how much they can contribute together.

\subsection{Following a proximal trajectory on one time interval}
\label{app:tracking}

Let $F$ agree with $f$ on $\cX$ and equal $+\infty$ elsewhere. For $\lambda>0$ and $x\in\cX$,
let $P_\lambda x$ be the unique minimizer over $z\in\cX$ of $f(z)+\|z-x\|^2/(2\lambda)$, and let
\[
 R_\lambda(x)\coloneqq\|x-P_\lambda x\|,\qquad
 M_\lambda(x)\coloneqq f(P_\lambda x)+\frac{R_\lambda(x)^2}{2\lambda}.
\]
The displacement $R_\lambda(x)$ will control changes of the minimum value $M_\lambda(x)$ as $\lambda$ varies.
By proximal optimality, $(x-P_\lambda x)/\lambda\in\partial F(P_\lambda x)$. Applying monotonicity
of $\partial F$ at two proximal points gives firm nonexpansiveness of $P_\lambda$, so both $P_\lambda$
and $x\mapsto x-P_\lambda x$ are $1$-Lipschitz. In particular, $R_\lambda$ is $1$-Lipschitz and
vanishes on $\cX_\star$. If $v\in\partial F(x)$, monotonicity at $x$ and $P_\lambda x$ also gives
$R_\lambda(x)^2\le\lambda\langle v,x-P_\lambda x\rangle$. By the bounded selection constructed in
Appendix~\ref{gd:app:upper}, it follows that $R_\lambda(x)\le\lambda$ for every $x\in\cX$.

Fix integers $s,k,m$ with $1\le s\le n$, $k\ge1$, and $0\le m\le\min\{k,n-s\}$. Starting from
$z_0\coloneqq X_s$, let $z_{j+1}\coloneqq P_1z_j$ for $0\le j<m$ and let
$H\coloneqq\max_{0\le j\le m}\|X_{s+j}-z_j\|$. We first show that $\|H\|_{2,s}\lesssim\sqrt k$.
If the innovations on this interval have norm at most $2\sqrt k$, we also show that $\|H\|_{4,s}\lesssim\sqrt k$.
For $0\le j\le m$ let $a_j\coloneqq X_{s+j}-z_j$, and for $0\le j<m$ let $v_j\coloneqq z_j-z_{j+1}$ and
$u_j\coloneqq a_j-G_{s+j}+v_j$. Since $v_j\in\partial F(z_{j+1})$, monotonicity gives
$\langle a_j+v_j,G_{s+j}-v_j\rangle\ge0$. Expanding the square and using $\|G_{s+j}\|\le1$ yields
$\|u_j\|^2\le\|a_j\|^2+\|G_{s+j}\|^2-\|v_j\|^2\le\|a_j\|^2+1$. Projection therefore gives
\[
 \|a_{j+1}\|^2\le\|a_j\|^2+1+\|\xi_{s+j}\|^2+2\langle u_j,\xi_{s+j}\rangle .
\]
Since $u_j$ is measurable before the innovation, conditioning and summing gives
$\E_s[\|a_j\|^2]\le2j$. For $0\le l\le m$ let $M_l'\coloneqq\sum_{j=0}^{l-1}\langle u_j,\xi_{s+j}\rangle$.
Orthogonality of its martingale increments gives $\E_s[|M_m'|^2]\le\sum_{j=0}^{m-1}(2j+1)\le k^2$.
By Doob's $L^2$ maximal inequality the conditional second-moment norm of
$\max_{0\le l\le m}|M_l'|$ is at most $2k$. Summing the pathwise inequality and using
$\E_s[\sum_{j=0}^{m-1}\|\xi_{s+j}\|^2]\le k$ proves the first estimate.

Under the additional bound $\|\xi_{s+j}\|\le2\sqrt k$, write $w_j\coloneqq\|\xi_{s+j}\|^2$.
Since $w_j^2\le4k w_j$ and $\E_{s+j}[w_j]\le1$, the diagonal terms of
$\E_s[(\sum_{j=0}^{m-1}w_j)^2]$ sum to at most $4km$, and conditioning at the later index bounds
each off-diagonal term by one. Thus $\|\sum_{j=0}^{m-1}w_j\|_{2,s}\le\sqrt5\,k$.
Taking conditional second-moment norms of the pathwise bound for $H^2$ now proves the fourth-moment
estimate. The same estimates hold for the maximum of partial sums of centered innovations multiplied
by predictable scalar weights of absolute value at most one, by taking the zero objective on $\R^d$.

We now remove the additional bound on the innovations. Fix one update $t$ in the interval and let
$b\coloneqq\sqrt k$, $V\coloneqq\xi_t\mathbf1_{\{\|\xi_t\|\le b\}}$,
$\mu\coloneqq\E_t[V]$, and $e_t\coloneqq\E_t[\|\xi_t\|^2\mathbf1_{\{\|\xi_t\|>b\}}]$.
The moment bound gives $e_t\le v\coloneqq k^{-\varepsilon/2}$, and centering gives
$\|\mu\|\le e_t/b$. Truncating at $b$ leaves mean $\mu$. To restore centering while changing the
innovation on as few paths as possible, take a fresh independent uniform variable $U$ on $[0,1]$
and let
\[
 \widehat\xi_t\coloneqq V-b\frac{\mu}{\|\mu\|}
       \mathbf1_{\{U\le\|\mu\|/b\}},\qquad \delta_t\coloneqq\xi_t-\widehat\xi_t ,
\]
taking the correction to be zero when $\mu=0$. We use a separate uniform for each comparison and
update, independent of the original process, and reveal it with that update's innovation. All the
original conditional bounds therefore persist in the enlarged filtration. The correction has
conditional mean $-\mu$, conditional second moment $b\|\mu\|$, and is conditionally independent
of $\xi_t$. Consequently $\widehat\xi_t$ and $\delta_t$ are centered,
$\|\widehat\xi_t\|\le2\sqrt k$, and
\[
 \E_t[\|\widehat\xi_t\|^2]\le1,\qquad
 \E_t[\|\delta_t\|^2]\le2v,\qquad \Pp_t[\delta_t\ne0]\le2v/k .
\]
For the first bound the squared norm has conditional expectation
$\E_t[\|V\|^2]+b\|\mu\|-2\|\mu\|^2\le\E_t[\|\xi_t\|^2]$. For the second, replace the first term
by $e_t$, and for the last bound add the probabilities of truncating and correcting.

Run a second projected trajectory from $X_s$ using $\widehat\xi_t$, choosing the original
subgradient whenever the two points agree and the bounded measurable selection from
Appendix~\ref{gd:app:upper} elsewhere. For $0\le j\le m$ let $d_j$ be the difference of the two
points at time $s+j$, and for $j<m$ let $h_j$ be the corresponding difference of subgradients.
Until the first change in the noise the points agree, so the conditional probability given $\cF_s$
that $h_j\ne0$ is at most $2jv/k$. Since $\|h_j\|\le2$ and $\langle d_j,h_j\rangle\ge0$,
nonexpansiveness of projection gives
\[
 \|d_{j+1}\|^2\le\|d_j\|^2+\|h_j\|^2+\|\delta_{s+j}\|^2
                   +2\langle d_j-h_j,\delta_{s+j}\rangle .
\]
After conditioning and summing, $\E_s[\|d_j\|^2]\lesssim v(j+j^2/k)\lesssim vj$. Moreover,
$\|d_j-h_j\|^2\le\|d_j\|^2+\|h_j\|^2$, so orthogonality bounds the conditional second moment
of the sum of the last inner products by $2v\sum_{j=0}^{m-1}\E_s[\|d_j-h_j\|^2]\lesssim k^2v^2$.
Applying Doob's inequality to that sum and summing the two squared terms yields
$\E_s[\max_{0\le j\le m}\|d_j\|^2]\lesssim kv$.

Both trajectories start at $X_s$, so we compare them to the same proximal orbit. By the triangle
inequality, the tracking error of the original path is at most the tracking error of the truncated
path plus the maximum distance between the two paths. Thus, for nonnegative random variables
$U_k,V_k$ on the enlarged space, we have
\begin{equation}\label{eq:tracking-pieces}
 H/\sqrt k\le U_k+V_k,\qquad
 \|U_k\|_{4,s}\lesssim1,\qquad \|V_k\|_{2,s}\lesssim k^{-\varepsilon/4}.
\end{equation}
The separate uniforms allow us to use these inequalities on many intervals of the same trajectory.
All subsequent estimates use moment inequalities, so the variables for overlapping intervals may
depend on one another.

\subsection{Controlling the motion near the last iterate}
\label{app:prox-motion}

We next use tracking to bound $R_N(X_t)$, the distance from the current iterate to its proximal
point at scale $N$. We will need a bound that also records how much of the initial distance can
still affect the last few updates.

First fix $x\in\cX$, let $z_0\coloneqq x$, and for every integer $j\ge0$ let $z_{j+1}\coloneqq P_1z_j$. If $y\in\cX$ has a
subgradient of norm at most $b$, the first proximal step from $y$ has length at most $b$.
By nonexpansiveness of $P_1$, successive steps have nonincreasing lengths and the two orbits started
from $x$ and $y$ stay within distance $\rho\coloneqq\|x-y\|$. It follows that
$\|z_j-y\|\le\rho+jb$ and $\|z_j-x\|\le2\rho+jb$.
For an integer $k\ge2$, now take $y=P_kx$, so $b=\rho/k$. We claim that
$\|z_k-z_{k-1}\|\le3\rho/k$. Let $l\coloneqq\lfloor k/2\rfloor$. For $1\le j\le l$, proximal
optimality and the subgradient inequality give
\[
 2(f(z_j)-f(y))\le\|z_{j-1}-y\|^2-\|z_j-y\|^2-\|z_j-z_{j-1}\|^2 .
\]
Since $f(z_j)$ decreases, summing gives $f(z_l)\le f(y)+\rho^2/(2l)$. The subgradient at $y$
gives $f(z_k)\ge f(y)-b(\rho+kb)$. Each step decreases $f$ by at least its squared length,
as follows by testing the subgradient $z_{j-1}-z_j$ at $z_{j-1}$. Thus the last $k-l$ steps give
$(k-l)\|z_k-z_{k-1}\|^2\le\rho^2/(2l)+2\rho^2/k$. Since $l\ge k/3$ and $k-l\ge k/2$,
this proves the claim.

Fix integers $N\ge1$ and $1\le t\le n$. To control $R_N(X_t)/\sqrt N$, we look back over successively
longer windows. For a window of $k$ updates with $2\le k\le t-1$, the last step of the proximal orbit
started at $X_{t-k}$ supplies a subgradient of norm at most $3R_k(X_{t-k})/k$ at its endpoint. Since
$R_N$ there is at most $N$ times this norm, the coefficient of the earlier normalized displacement
$R_k(X_{t-k})/\sqrt k$ is $3\sqrt{N/k}$. We choose $k=36N$ whenever the window fits, so each time we
enlarge the scale by a factor of $36$, we halve the contribution of the earlier displacement. Since
$1/2\le36^{-1/6}$, we will bound the contribution of the initial distance by a universal multiple of
\[
 I_{N,t}\coloneqq\frac{D}{\sqrt{t-1+N}}\left(\frac{N}{t-1+N}\right)^{1/6}.
\]
Here and below we allow the nonnegative random variables $U$ and $V$ to depend on the quantity and
interval being estimated. We claim that for a universal $C$ we can choose them so that
\begin{equation}\label{eq:prox-split}
 R_N(X_t)/\sqrt N\le U+V+C I_{N,t},\qquad
 \|U\|_4\le C,\qquad \|V\|_2\le C N^{-\varepsilon/4}.
\end{equation}
When $36N\le t-1$, let $H$ be the tracking error over the $36N$ updates from $X_{t-36N}$ to $X_t$.
Applying the bound at the endpoint of the proximal orbit and using the $1$-Lipschitz property of $R_N$
to return to $X_t$ gives
\[
 \frac{R_N(X_t)}{\sqrt N}
 \le\frac12\,\frac{R_{36N}(X_{t-36N})}{\sqrt{36N}}
       +6\,\frac{H}{\sqrt{36N}} .
\]
Apply this inequality recursively. Starting with $N_0=N$ and $t_0=t$, let
$N_{j+1}\coloneqq36N_j$ and $t_{j+1}\coloneqq t_j-36N_j$ until the first index $h$ with
$t_h-1<36N_h$. On this last interval start the proximal orbit at $X_1$ and follow it to time $t_h$,
using the window bound $36N_h$ in \eqref{eq:tracking-pieces}. Since the nearest minimizer is a fixed
point of $P_1$, the orbit stays within distance $D$ of it. The remaining residual is therefore at most
$D+H$. The coefficients of the normalized tracking errors are bounded by $6\cdot2^{-j}$, whose
sum is at most $12$. By Minkowski's inequality, \eqref{eq:tracking-pieces} gives the stated norms of
$U,V$, because every window has length bound at least $N$. Finally, since $N_h=36^h N$, we have $2^{-h}\le(N/N_h)^{1/6}$. Summing the elapsed
windows gives $N_h\le t-1+N\le39N_h$, hence $2^{-h}D/\sqrt{N_h}\lesssim I_{N,t}$, which proves
\eqref{eq:prox-split}.

We now specialize the bound to a suffix of $N$ updates. For $1\le N\le n-1$ let
$s\coloneqq n-N$, $A\coloneqq D^2/n$, $a_N\coloneqq\sqrt A(N/n)^{1/6}$, and
$Z_N\coloneqq\max_{s\le t\le n}\|X_t-X_s\|$. The orbit from $X_s$ moves at most $3R_N(X_s)$
in $N$ steps by the deterministic motion estimate with $y=P_N X_s$. Since $s-1+N=n-1\ge n/2$,
\eqref{eq:tracking-pieces} and \eqref{eq:prox-split} show that $Z_N/\sqrt N$
is at most $U+V+C a_N$ with the same moment bounds as in \eqref{eq:prox-split}, up to universal
constants. For every integer $N/2\le l\le N$, we also have $s-1+l\ge(n-1)/2$.
In particular,
\begin{equation}\label{eq:geom-second}
 \E[Z_N^2]\lesssim N(1+a_N^2),\qquad
 \E[R_l(X_s)^2]\lesssim l(1+a_N^2).
\end{equation}

Two further uses of these estimates will be convenient. First, for integers $1\le\lambda\le2n$,
\eqref{eq:prox-split} at time $n$ bounds $R_\lambda(X_n)/\sqrt\lambda$ by $U+V+C a_\lambda$,
where $a_\lambda\coloneqq\sqrt A(\lambda/n)^{1/6}$, $\|U\|_4\lesssim1$ and
$\|V\|_2\lesssim\lambda^{-\varepsilon/4}$. Second, for any predictable scalar weights
$(w_t)_{t=s}^{n-1}$ with $|w_t|\le1$, let $m\coloneqq N^{-1/2}\sum_{t=s}^{n-1} w_t\xi_t$.
By orthogonality $\E_s[\|m\|^2]\le1$. Replacing its innovations by their centered modifications at
cutoff $\sqrt N$ gives a sum with bounded fourth-moment norm, by the zero-objective case of the
tracking estimate. The difference has second-moment norm at most $\sqrt2\,N^{-\varepsilon/4}$,
again by orthogonality. Hence $\|m\|\le U+V$ with the same two moment bounds.

\subsection{Comparing two intervals ending at the last iterate}
\label{app:suffixes}

Let $x\coloneqq X_n$. To compare different scales, we need two pieces of information from a suffix of
$N$ updates. The average of its early iterates gives a point with objective value close to a
smoothed value at $x$, while the drift over the full suffix supplies a supporting inequality that
we can test at a point from an earlier interval.

Fix $1\le N\le n-1$ and an integer $N/2\le l\le N$. Let $s\coloneqq n-N$, and define the average
$\bar X\coloneqq l^{-1}\sum_{t=s}^{s+l-1}X_t$. For $s\le t<n$ let
$w_t\coloneqq\min\{1,(t-s+1)/l\}$ and let $\nu\coloneqq\sum_{t=s}^{n-1}w_t=N-(l-1)/2$.
Averaging the identities $x-X_j=\sum_{t=j}^{n-1}(X_{t+1}-X_t)$ over $s\le j<s+l$ gives these
weights. Define the effective drift $b_t\coloneqq X_t-X_{t+1}+\xi_t$, its weighted average
$\Gamma\coloneqq\nu^{-1}\sum_{t=s}^{n-1}w_t b_t$, and the normalized vectors
\[
 v\coloneqq\frac{x-\bar X}{\sqrt N},\qquad
 m\coloneqq\frac1{\sqrt N}\sum_{t=s}^{n-1}w_t\xi_t,\qquad
 y\coloneqq m-v=\frac{\nu}{\sqrt N}\Gamma .
\]
Because $\|v\|\le2Z_N/\sqrt N$, the estimates of Appendix~\ref{app:prox-motion} imply that
each of $\|v\|,\|y\|$ is bounded by $U+V+C a_N$, with $\|U\|_4\lesssim1$ and
$\|V\|_2\lesssim N^{-\varepsilon/4}$. In particular,
$\E[\|v\|^2]+\E[\|y\|^2]\lesssim1+a_N^2$.

Let $\lambda\coloneqq N/8$ and $h\coloneqq M_\lambda(x)$. We will construct nonnegative errors
$C_s,E_s$ with $\E[C_s+E_s]\lesssim1+a_N^2$ such that, for every $z\in\cX$,
\begin{equation}\label{eq:geom-support}
 f(\bar X)\le h+C_s,\qquad f(z)\ge h+\langle\Gamma,z-x\rangle-E_s .
\end{equation}
For the first inequality take $z_s\coloneqq P_lX_s$, which is measurable before any update of the
suffix. Summing the projected squared-distance inequality with comparator $z_s$ and applying
convexity gives
\[
 f(\bar X)\le M_l(X_s)+\frac1{2l}\sum_{t=s}^{s+l-1}\|G_t-\xi_t\|^2+
              \frac1l\sum_{t=s}^{s+l-1}\langle\xi_t,X_t-z_s\rangle .
\]
For every $z\in\cX$, the squared triangle inequality with weights $\lambda$ and $l-\lambda$
gives $\|z-X_s\|^2/l\le\|z-x\|^2/\lambda+\|x-X_s\|^2/(l-\lambda)$.
Testing $P_\lambda x$ in the definition of $M_l(X_s)$ therefore gives
$M_l(X_s)\le h+Z_N^2/(2(l-\lambda))$. Take $C_s$ to be the latter error plus the last two terms of
the display, with the martingale sum in absolute value. By \eqref{eq:geom-second}, the first error
has expectation at most a constant times $1+a_N^2$, and conditional centering bounds the expected squared
update norm by two. Orthogonality bounds the second moment of the martingale sum by
$2l\,\E[Z_N^2+R_l(X_s)^2]$. Taking the square root and dividing by $l\ge N/2$ bounds the expected
absolute value of the normalized sum by a constant times $\sqrt{1+a_N^2}$, which suffices.

For the second inequality, projection gives
$\langle b_t-G_t,z-X_{t+1}\rangle\le0$ for $z\in\cX$, and both the projection residual
$\|b_t-G_t\|$ and the step length $\|X_{t+1}-X_t\|$ are at most $1+\|\xi_t\|$.
The subgradient inequality at $X_t$ therefore gives
$f(z)\ge f(X_t)+\langle b_t,z-X_t\rangle-(1+\|\xi_t\|)^2$.
To bound the objective term from below, let $g\coloneqq(x-P_\lambda x)/\lambda$. Testing this
subgradient at $X_t$ gives $f(X_t)\ge h+\langle g,X_t-x\rangle+\lambda\|g\|^2/2$.
Completing the square and using $\|X_t-x\|\le2Z_N$ gives $f(X_t)\ge h-2Z_N^2/\lambda$.
Averaging the resulting supporting inequalities with weights $w_t/\nu$ now gives
\eqref{eq:geom-support} with
\[
 E_s\coloneqq\frac{2Z_N^2}{\lambda}+\|\Gamma\|Z_N+
 \frac1\nu\sum_{t=s}^{n-1}w_t(1+\|\xi_t\|)^2+
 \frac1\nu\left|\sum_{t=s}^{n-1}w_t\langle b_t,X_t-X_s\rangle\right|.
\]
The first three terms have total expectation of order at most $1+a_N^2$ by
\eqref{eq:geom-second}, $\nu\ge N/2$, and the second-moment bound for $y$.
In the last term write $b_t=\xi_t-(X_{t+1}-X_t)$. The martingale sum has second moment at
most $N\E[Z_N^2]$. For the remaining sum use
$2\langle X_{t+1}-X_t,X_t-X_s\rangle
=\|X_{t+1}-X_s\|^2-\|X_t-X_s\|^2-\|X_{t+1}-X_t\|^2$.
Since the weights increase from a nonnegative value to one, summation by parts bounds the absolute
value of this remaining sum by $Z_N^2+\tfrac12\sum_{t=s}^{n-1}(1+\|\xi_t\|)^2$.
Dividing by $\nu$ proves the required estimate for $\E[E_s]$.

We now apply \eqref{eq:geom-support} on nested suffixes. Let $r\ge2$ be a power of two such that
$P\coloneqq\lfloor\log(n-1)/(2\log r)\rfloor\ge1$. For $0\le i\le P$ let
$N_i\coloneqq r^{2(P-i)}$, $s_i\coloneqq n-N_i$, and $h_i\coloneqq M_{N_i/8}(x)$.
For $i<P$ take $l_i\coloneqq N_i-N_{i+1}$, and take $l_P\coloneqq1$. The first $l_i$ points in
suffix $i$ precede the next suffix, so their average is known before its noise is drawn. Use the
subscript $i$ for $\bar X,v,m,y,\nu,\Gamma,C_s,E_s$ from that choice. For $0\le i<j\le P$,
applying the supporting inequality from suffix $j$ at $\bar X_i$ and the upper inequality from
suffix $i$ gives
\begin{equation}\label{eq:geom-cross}
 \langle v_i,y_j\rangle\ge\frac{\nu_j}{N_j}r^{i-j}
                     (h_j-h_i-C_i-E_j).
\end{equation}

To bound the rises of the envelope between consecutive suffixes, the determinant argument in
Appendix~\ref{app:determinant} will collect earlier displacement vectors in a single matrix. Expanding
its quadratic form pairs a later drift with every displacement already recorded, so we now bound the
sum of $|\langle v_i,y_j\rangle|$ over all $0\le i<j\le P$. Since decreasing the proximal parameter increases the envelope,
$h_j\ge h_i$. Also $\nu_j/N_j\le1$ and the weights $r^{i-j}$ sum to at most one over either index
when the other is fixed. Thus \eqref{eq:geom-cross} bounds the expected sum of negative parts by
$\sum_{i=0}^P\E[C_i+E_i]\lesssim P+A$, using
$\sum_{i=0}^P a_{N_i}^2\lesssim A$ from the geometric progression of lengths.
For the signed sum, write $y_j=m_j-v_j$. The average $\bar X_i$ and the point $X_{s_j}$ are both
$\cF_{s_j}$-measurable, and $\E_{s_j}[m_j]=0$. Hence
\[
 |\E[\langle v_i,m_j\rangle]|
 =\frac{|\E[\langle x-X_{s_j},m_j\rangle]|}{\sqrt{N_i}}
 \lesssim r^{i-j}(1+a_{N_j}).
\]
Here we used $\E[\|m_j\|^2]\le1$ and \eqref{eq:geom-second}. Summing this bound costs at most
a constant times $P+\sqrt A$, since $\sum_{j=0}^P a_{N_j}\lesssim\sqrt A$.
Expanding $\|\sum_{i=0}^P v_i\|^2\ge0$ bounds $-\sum_{0\le i<j\le P}\langle v_i,v_j\rangle$
by $\tfrac12\sum_{i=0}^P\|v_i\|^2$. Since $P\ge1$, the signed sum therefore has expectation
at most a constant times $P+A$. Adding twice the negative parts proves
\begin{equation}\label{eq:geom-triangle}
 \E\!\left[\sum_{0\le i<j\le P}|\langle v_i,y_j\rangle|\right]\lesssim P+A .
\end{equation}

\subsection{Summing the large noise increments across scales}
\label{app:scale-tails}

It remains to bound the rise of the envelope from the longest suffix to the last update.
For $0\le i<P$ let $\delta_i\coloneqq h_{i+1}-h_i$. A second use of the proximal displacement
will let us pay for unusually large increments at the scales where they occur. For each $\lambda>0$,
testing $P_\lambda x$ at parameter $\lambda/2$ gives
\[
 0\le M_{\lambda/2}(x)-M_\lambda(x)\le R_\lambda(x)^2/(2\lambda).
\]
Let $\Lambda\coloneqq2^{\lceil\log_2(n-1)\rceil}$, where $\log_2$ is the base-two logarithm.
Starting at $\Lambda$ and halving down to $1/8$ includes all parameters $N_i/8$.
For every parameter $\lambda\ge1$ on this grid, squaring the terminal estimate of
Appendix~\ref{app:prox-motion} bounds the increment by two nonnegative quantities
$\mathcal U_\lambda,\mathcal V_\lambda$ satisfying
\[
 \|\mathcal U_\lambda\|_2\lesssim1,\qquad
 \E[\mathcal V_\lambda]\lesssim\lambda^{-\varepsilon/2}+A(\lambda/n)^{1/3}.
\]
In fact, we may use a universal multiple of $U^2$ for the first term and of
$V^2+a_\lambda^2$ for the second. For $\lambda<1$, $R_\lambda(x)\le\lambda$ gives the deterministic
bound $\lambda/2$, which we assign to $\mathcal U_\lambda$.

The sum of the expected $\mathcal V_\lambda$ over the whole grid is of order at most
$\Phi_{n,2+\varepsilon}+A$. To see why this is uniform down to $\varepsilon=0$, let
$h\coloneqq\log_2\Lambda$. Comparing the decreasing summand to its integral gives
\[
 \sum_{j=0}^{h}2^{-j\varepsilon/2}
 \le1+\frac{2}{\log 2}\int_1^{\sqrt\Lambda}t^{-1-\varepsilon}\,dt
 \lesssim\Phi_{n,2+\varepsilon}.
\]
Here $\sqrt\Lambda\le n$. The sum of $(\lambda/n)^{1/3}$ is bounded by a numerical constant
by geometric summation from $\Lambda\le2n$ downward. For each coarse interval from $N_i/8$ to
$N_{i+1}/8$, denote by $U_i,V_i$ the sums of the corresponding fine bounds. Because there are
$2\log_2 r$ increments in that interval, Minkowski's inequality gives $\|U_i\|_2\lesssim\log r$.
The initial range from $\Lambda$ to $N_0/8$ contains at most $4+2\log_2 r$ increments, since
$N_0\le n-1<r^2N_0$. Its contributions from $\mathcal U_\lambda$ thus have total expectation
of order at most $\log r$.

To start the sum, choose a minimizer $x_\star$ nearest $X_1$. Since
$\langle G_t,X_t-x_\star\rangle\ge0$, expanding the projected update and summing gives
$\E[\|x-x_\star\|^2]\le D^2+2(n-1)$. Testing $x_\star$ gives
$\E[M_\Lambda(x)]\lesssim A+1$. At the other end, Lipschitz continuity gives
$f(x)-M_{1/8}(x)\le1/16$ by completing the square. For each $0\le i<P$ let
$\delta_i'\coloneqq\min\{\delta_i,U_i\}$. Since $\delta_i\le U_i+V_i$, we have
$\delta_i\le\delta_i'+V_i$ and $\|\delta_i'\|_2\lesssim\log r$. Summing over the disjoint ranges
of the grid now gives
\begin{equation}\label{eq:geom-grid}
 \E[f(x)]\lesssim A+\Phi_{n,2+\varepsilon}+\log r+
                     \sum_{i=0}^{P-1}\E[\delta_i'] .
\end{equation}

We will apply a dimension-dependent estimate to pairs of vectors with norms at most $r$.
For $0\le i<P$, let $B_i\coloneqq\{\|v_i\|>r\text{ or }\|y_{i+1}\|>r\}$. Each norm is bounded
by $U+V+C a_N$ as in Appendix~\ref{app:suffixes}. Applying Markov's inequality separately to the
three terms at threshold $r/3$, using the fourth moment of $U$ and second moment of $V$, gives
\[
 \sqrt{\Pp[B_i]}\lesssim
 \frac1{r^2}+\frac{N_{i+1}^{-\varepsilon/4}}r+\frac{a_{N_i}}r .
\]
The deterministic term follows from $\mathbf1_{\{C a_N>r/3\}}\le 3C a_N/r$. Cauchy-Schwarz
therefore bounds $\E[\delta_i'\mathbf1_{B_i}]$ by the right-hand side times a universal multiple
of $\log r$. The first terms sum to order at most $P$, and the third to order at most $\sqrt A$
because $(\log r)/r\le1$. For the second terms we use the actual scales in the sum. Since
$r^{P/2}\le n$, comparison to an integral gives, also at $\varepsilon=0$,
\[
 \frac{\log r}{r}\sum_{i=0}^{P-1}N_{i+1}^{-\varepsilon/4}
 =\frac{\log r}{r}\sum_{j=0}^{P-1}r^{-j\varepsilon/2}
 \le\frac{\log r}{r}+\frac2r\int_1^{r^{P/2}}t^{-1-\varepsilon}\,dt
 \lesssim\Phi_{n,2+\varepsilon}.
\]
We have thus bounded the contribution of all events $B_i$ to the last sum in
\eqref{eq:geom-grid} by a constant times $P+A+\Phi_{n,2+\varepsilon}$.

\subsection{Where the dimension enters}
\label{app:determinant}

We now use \eqref{eq:geom-cross} on consecutive intervals. Let
$\widetilde v_i\coloneqq v_i\mathbf1_{\{\|v_i\|\le r\}}$ and
$\widetilde y_i\coloneqq y_i\mathbf1_{\{\|y_i\|\le r\}}$ for $0\le i\le P$.
Since $\nu_{i+1}/N_{i+1}\ge1/2$, on $B_i^c$ we have
$\delta_i'\le C_i+E_{i+1}+2r|\langle\widetilde v_i,\widetilde y_{i+1}\rangle|$.
The expected sum of $C_i+E_{i+1}$ is of order at most $P+A$. To bound the remaining inner
products we keep track of the directions already used.

Let $I_d$ be the $d$-dimensional identity matrix. For $0\le i\le P$ let
$A_i\coloneqq r^2 I_d+\sum_{k=0}^{i-1}\widetilde v_k\widetilde v_k^\top$.
For each $0\le i<P$, expanding $\widetilde y_{i+1}^{\top} A_i\widetilde y_{i+1}$ gives $r^2\|\widetilde y_{i+1}\|^2$
and one squared inner product for each $0\le k<i$. Since both truncated vectors have norm at most $r$, each such term satisfies
$\langle\widetilde v_k,\widetilde y_{i+1}\rangle^2
 \le r^2|\langle v_k,y_{i+1}\rangle|$. Each ordered pair occurs at most once in the sum.
Thus \eqref{eq:geom-triangle} and the second-moment bound for $y_{i+1}$ give
\[
 \E\!\left[\sum_{i=0}^{P-1}\widetilde y_{i+1}^{\top}A_i\widetilde y_{i+1}\right]
 \lesssim r^2(P+A).
\]
The other factor measures the length of $\widetilde v_i$ in the inverse matrix. For $0\le i<P$ let
$\theta_i\coloneqq\widetilde v_i^\top A_i^{-1}\widetilde v_i$. Since $A_i\succeq r^2I_d$,
meaning that $A_i-r^2I_d$ is positive semidefinite, $0\le\theta_i\le1$. The rank-one determinant
identity gives $\det A_{i+1}=(1+\theta_i)\det A_i$. Since $\theta_i\le2\log(1+\theta_i)$ and
$\operatorname{tr}A_P\le r^2(d+P)$, the arithmetic-geometric mean inequality for the eigenvalues yields
\[
 \sum_{i=0}^{P-1}\theta_i
 \le2\log\frac{\det A_P}{\det A_0}\le2d\log(1+P/d).
\]
Cauchy-Schwarz in the two matrix norms and then in the sum over $i$ bounds the sum of adjacent
absolute inner products by the square root of the product of the two sums above. The determinant
bound holds on every path, so taking expectations gives an upper bound of order
$r\sqrt{(P+A)d\log(1+P/d)}$. Combining with \eqref{eq:geom-grid} and the bound on exceptional
events gives
\begin{equation}\label{eq:geom-budget}
 \E[f(X_n)]\lesssim A+\Phi_{n,2+\varepsilon}+P+\log r+
                  r^2\sqrt{(P+A)d\log(1+P/d)} .
\end{equation}

It remains to choose the ratio of successive scales. Let $J\coloneqq\log n$ and $w\coloneqq J/d$.
If $w\le2^{16}$, then $\log(e+w)$ is bounded by a numerical constant. Applying the variance estimate
\eqref{gd:eq:variance-normalized} gives $\E[f(X_n)]\lesssim A+1+J\lesssim A+1+q_{n,d}$.
Suppose now that $w>2^{16}$ and choose $r\coloneqq2^{\lfloor(\log_2 w)/8\rfloor}$.
Then $r\ge2$, $\log r\asymp\log w$, and $\log(n-1)\asymp J$. In particular
$\log(n-1)/(2\log r)\ge2$, so taking its integer part gives
$P\asymp J/\log w\asymp q_{n,d}$. Since $(\log w)^2\lesssim w\le J$, we also have $\log r\lesssim P$.
Finally $r^4\le w^{1/2}$, $d/P\lesssim(\log w)/w$, and $P/d\lesssim w$, whence
\[
 \frac{r^4 d\log(1+P/d)}{P}\lesssim
 \frac{(\log w)^2}{\sqrt w}\lesssim1.
\]
The last term of \eqref{eq:geom-budget} is therefore at most a constant times
$\sqrt{(P+A)P}$, which is at most $P+A$. In either range of $w$, we obtain
$\E[f(X_n)]\lesssim D^2/n+\Phi_{n,2+\varepsilon}+q_{n,d}$.
For $2\le p\le3$ the moment term is exactly $\Phi_{n,p}$, and for $3<p\le\infty$ we have
$\Phi_{n,3}\le2\le2\Phi_{n,p}$. All comparison constants are thus independent of $p$.

For general $\eta>0$, apply the unit-step result on $\cX/\eta$ to the points $X_k/\eta$ and the
objective $z\mapsto(f(\eta z)-f_\star)/\eta$. Its Lipschitz and noise bounds are unchanged,
and the initial distance is at most $D/\eta$. Multiplying the conclusion by $\eta$ bounds the
original error by a constant times $D^2/(\eta n)+\eta b_{n,p,d}$. When $D\ge\eta n$, the bound on
each expected step length gives the simpler estimate
$\E[\Delta(X_n)]\le D+2\eta(n-1)\le3D$. Taking the applicable bound on either side of $D=\eta n$
proves the upper estimate in \eqref{eq:main-prescribed}.

\section{Lower bounds from scalar instances}
\label{sec:lower}

We use two families of linear objectives on the half-line $[0,\infty)$, one to measure how far a noise
increment can displace the last iterate, the other to measure the effect of the initial distance. In
dimension $d>1$, use the product domain $[0,\infty)\times\R^{d-1}$, let the objective depend only on the
first coordinate, and append zero coordinates to the starting point, selected subgradients and noises.
Projection then gives exactly the scalar recursion in the first coordinate.

\subsection{Independent noise starting at a minimizer}
\label{sec:lower-noise}

Fix $n\ge2$, $p\in[1,\infty]$ and $\sigma>0$. We construct a linear objective of slope at most $\sigma$,
started at $X_1=0$, such that for every $\eta>0$,
\begin{equation}
 \E[\Delta(X_n)]\ge\frac{\eta\sigma^2}{128}\Phi_{n,p}.
 \label{eq:lower-noise}
\end{equation}
Following the time-dependent construction of Beretta et al.~\cite{BCCP2026}, we choose the size of a
negative noise increment according to how many updates remain. With $j$ updates left, a jump of order
$\eta\sigma j$ in the iterate can survive the remaining steps back toward zero. Its allowed probability
is governed by the $p$-th moment, which is why the eventual sum involves $j^{1-p}$. The noises are independent and the
filtration is generated by the preceding noises, so centering and moment bounds for the individual laws
give the conditional assumptions.

\emph{Moment orders $1\le p\le2$.}
Let $f(x)\coloneqq\sigma x$ and use the subgradient $\sigma$ at every update. When $p$ is close to one, the sum of the candidate jump probabilities can be large. We normalize these
probabilities so that, conditional on any particular jump, there is a uniformly positive chance it is
the only one. Let $K\coloneqq(\sum_{j=1}^{n-1} j^{-p})^{1/p}$. For $1\le j<n$, set $A_j\coloneqq4\sigma Kj$,
$r_j\coloneqq1/(2(4Kj)^p)$ and $b_j\coloneqq r_j A_j/(1-r_j)$. At time $n-j$, let the noise take value
$-A_j$ with probability $r_j$ and value $b_j$ otherwise. The chosen $b_j$ centers the noise, and since
$K\ge1$ and $r_j\le1/8$,
\[
 \E[|W_{n-j}|^p]=r_j A_j^p\left[1+\left(\frac{r_j}{1-r_j}\right)^{p-1}\right]
 \le2r_j A_j^p=\sigma^p .
\]
Let $E_j$ be the event that only the noise at time $n-j$ takes its negative value. Since $\sum_{i=1}^{n-1}
r_i=1/(2\cdot4^p)\le1/8$, independence and the union bound give $\Pp[E_j]\ge(7/8)r_j$. For $1\le i<n$
we have $b_i\le(8/7)r_i A_i\le4\sigma/7$. On $E_j$, the process stays at zero until update $n-j$. The
total of the $j$ deterministic subgradients and the remaining positive noises is at most
$(11/7)j\sigma\le A_j/2$, giving
\[
 X_n\ge\eta\left(A_j-j\sigma-\sum_{i=1}^{j-1}b_i\right)_+\ge\eta A_j/2 .
\]
The events are disjoint, hence $\E[\Delta(X_n)]\ge 7\eta\sigma^2 S_{n-1,p}/(32(4K)^{p-1})$. At $p=1$,
$K^{p-1}=1$. For $1<p\le2$, integral comparison gives $K^{p-1}\le(p/(p-1))^{(p-1)/p}\le e^{1/e}<2$. Thus
$(4K)^{p-1}\le8$, and \eqref{eq:sum-comparison} proves \eqref{eq:lower-noise}.

\emph{Moment orders $2<p<\infty$.}
Here $\sum_{j=2}^\infty j^{-p}$ is uniformly bounded, so we can use jump probabilities $1/(2j^p)$
directly. We halve the slope so that a negative jump at the very last update already creates a positive
error. Let $f(x)\coloneqq\sigma x/2$ and select its slope at each update. For $1\le j<n$, let
$r_j\coloneqq1/(2j^p)$ and choose $W_{n-j}$ to equal $-\sigma j$ with probability $r_j$ and $\sigma
j/(2j^p-1)$ otherwise. The mean is zero and the $p$-th absolute moment divided by $\sigma^p$ is
$1/2+1/(2(2j^p-1)^{p-1})\le1$. Now let $E_j$ denote the event that the first negative noise is at time
$n-j$. By integral comparison,
$\tfrac12\sum_{i=2}^\infty i^{-p}\le\tfrac12(1/4+1/9+1/16+1/25+\int_5^\infty t^{-2}\,dt)<1/3$, so
$\Pp[E_j]\ge1/(3j^p)$. On this event, the smallest possible last iterate is obtained by taking the
positive value at every subsequent update, which gives
\[
 X_n\ge\eta\sigma\left(\frac j2-\sum_{i=1}^{j-1}\frac{i}{2i^p-1}\right)_+ .
\]
For $j=1$ the bound is $\eta\sigma/2$. For $j\ge8$, the first seven summands total at most
$1+(4/7)\sum_{i=2}^7 i^{-1}<2$ and each subsequent summand is at most $1/4$, so the sum is at most $j/4$.
Using the disjoint events for $j=1$ and $8\le j<n$ gives
$\E[\Delta(X_n)]\ge\eta\sigma^2(1/12+\sum_{j=8}^{n-1} j^{1-p}/24)$. Since
$\sum_{j=2}^{\min\{7,n-1\}}j^{1-p}<2$, this is at least $\eta\sigma^2 S_{n-1,p}/64$. Applying
\eqref{eq:sum-comparison} proves \eqref{eq:lower-noise}.

\emph{Bounded noise.}
Use the same slope $\sigma/2$ and zero noise except at the last update, where the two signs of $\sigma$
are equiprobable. The final iterate is zero or $\eta\sigma/2$ with equal probability, so
$\E[\Delta(X_n)]=\eta\sigma^2/8$. Since $\Phi_{n,\infty}=1$, this proves \eqref{eq:lower-noise}. In all
three constructions the objective and noise laws depend only on $n,p,\sigma$, and projection onto the
half-line commutes with multiplication by $\eta$. Thus each constructed instance gives the claimed bound
simultaneously for all stepsizes.

\subsection{Initial distance and the resulting constants}

We now isolate the effect of the initial distance by choosing a slope small enough that even $n-1$
noiseless steps remove at most half that distance. For every $D>0$ and $\eta>0$, let
$a\coloneqq\min\{1/2,D/(2\eta n)\}$ and use the noiseless objective
$f(x)\coloneqq ax$ starting at $X_1=D$. Since $D-(n-1)\eta a\ge D/2$, selecting its slope keeps the
trajectory positive and gives $X_n=D-(n-1)\eta a$. The terminal error is therefore at least
$\tfrac14\min\{D,D^2/(\eta n)\}$. Combining this instance with \eqref{eq:lower-noise} at $\sigma=1$ and
embedding as above gives, in every dimension,
\begin{equation}
 \mathcal R_{n,p}^{(d)}(\eta,D)\ge
 \max\left\{\frac{\eta\Phi_{n,p}}{128},\
            \frac14\min\left\{D,\frac{D^2}{\eta n}\right\}\right\}.
 \label{eq:lower-all-steps}
\end{equation}
In particular, the sum of the distance term $\min\{D,D^2/(\eta n)\}$ and $\eta\Phi_{n,p}$ is at most
$132$ times the risk. If $\eta n<D$, the second lower bound in \eqref{eq:lower-all-steps} is $D/4$, and
$\Phi_{n,p}\le n$. Otherwise the geometric mean of the two lower bounds is
$(D/(16\sqrt2))\sqrt{\Phi_{n,p}/n}$. This proves the latter lower bound for the best constant stepsize as
well.

For $D>0$ on the line, \eqref{eq:normalized-upper} gives the matching upper bound $\sqrt{65}\,D\sqrt{\Phi_{n,p}/n}$
at $\eta=D/\sqrt{65n\Phi_{n,p}}$. At $D=1$ and $\eta=n^{-1/2}$, the lower and upper constants multiplying
$\Phi_{n,p}/\sqrt n$ are $1/128$ and $33$. In every dimension with $1\le p\le2$, using
\eqref{gd:eq:heavy-fixed-step} gives the optimized upper bound $4D\sqrt{\Phi_{n,p}/n}$ at
$\eta=D/(4\sqrt{n\Phi_{n,p}})$, and gives the standard-step upper constant $9$.

\section{A bounded-noise construction in the plane}
\label{gd:lb:appendix}

In this appendix $Z_t$ denotes the state after $t$ updates, so $Z_0$ is the starting point and $Z_{n-1}$ is the last iterate. We add the innovation at each update. When the
stepsizes and objectives agree, the additive innovation $\xi_{t+1}$ in the update from $Z_t$ corresponds
to $W_{t+1}=-\xi_{t+1}$ in \eqref{eq:ssgm}, with $X_{t+1}=Z_t$. Let $n_0\coloneqq1+65536\cdot10^6$ and, for
$n\ge n_0$, let $A_n\coloneqq\log(n-1)/\log\log(n-1)$.

\begin{lemma}\label{gd:thm:bounded-plane}
For every integer $n\ge n_0$ and stepsize $\rho>0$, there is a convex coercive function
$f:\R^2\to[0,\infty)$ with Lipschitz constant less than two, a unique minimizer at the origin, and
$f(0)=0$. It is a finite maximum of affine functions and admits a Borel selector $g:\R^2\to\R^2$ such
that $g(0)=0$ and, for every $z\in\R^2$, $g(z)\in\partial f(z)$. There are a filtration $(\cF_t)_{t=0}^{n-1}$
and innovations $(\xi_{t+1})_{t=0}^{n-2}$ such that for each $0\le t<n-1$, $\xi_{t+1}$ is
$\cF_{t+1}$-measurable, belongs to $\{0\}\times[-1,1]$ and satisfies $\E[\xi_{t+1}\mid\cF_t]=0$. Starting
at $Z_0=0$, the updates $Z_{t+1}\coloneqq Z_t-\rho g(Z_t)+\rho\xi_{t+1}$ satisfy
\begin{equation}
 \Pp[f(Z_{n-1})\ge\rho A_n/5760]\ge7/8,\qquad
 \E[f(Z_{n-1})]\ge\rho A_n/10000 .
 \label{gd:lb:plane-prob}
\end{equation}
\end{lemma}

We first work at unit stepsize. The difficulty is to raise the objective with centered noise when the selected vertical drift points
downward at positive heights. We use two convex profiles of the vertical coordinate. On each long
interval the horizontal coordinate makes the flatter profile supply the vertical subgradient, so the
noise can carry the iterate across against a small drift. The steeper profile gains a fixed amount on
that crossing. Short connecting intervals switch which profile supplies the subgradient before the next
crossing. We will bound the probability of completing all crossings and then choose how many fit in $n-1$
updates.

\subsection{The profiles and the subgradient selection}
\label{gd:lb:geometry}

The two profiles will take turns having the larger slope. We choose an odd number $m\ge1$ of long
intervals so that the same profile gains more on the first and last crossings. Let $M\ge1000$ be an
integer giving the length of the last interval. We choose the earlier intervals progressively longer and their slopes progressively
smaller, so that the difference between the two profile gains is the same at every crossing. Write $Q$ for the factor by which we shorten successive intervals. As we will see in the proof of
Lemma~\ref{gd:lb:stage-success}, conditional on all preceding stages succeeding, the probability that
a stage fails is at most a universal constant times $1/\sqrt Q$. To make the sum of these $m$ bounds
at most $1/8$, let $Q\coloneqq(256m)^2$. For $1\le j\le m$, let $\ell_j\coloneqq M Q^{m-j}$ be
the length of interval $j$.
Let $\delta_0\coloneqq1/8$. For $1\le j\le m$, define $b_j\coloneqq\delta_0/((Q-1)\ell_j)$, and let
$b_{m+1}\coloneqq Qb_m$ and $\beta\coloneqq b_{m+1}$. For every $1\le j\le m$ we then have $b_{j+1}=Qb_j$
and $(b_{j+1}-b_j)\ell_j=\delta_0$. In particular $\beta=Q\delta_0/((Q-1)M)\le2\delta_0/M$, so
$\beta<1/4$ and $4\beta<\delta_0/3$.

Place the intervals on the vertical axis with a connector of width $w\coloneqq10$ before each interval.
Let $a_1\coloneqq w$. For $1\le j\le m$, let $r_j\coloneqq a_j+\ell_j$, and if $j<m$, let
$a_{j+1}\coloneqq r_j+w$. For $1\le j\le m$, let $I_j\coloneqq[a_j,r_j]$ denote a long interval. The
first connector is $P_0\coloneqq[0,w]$, for $1\le j<m$ the intermediate connectors are
$P_j\coloneqq[r_j,a_{j+1}]$, and the final ray is $P_m\coloneqq[r_m,\infty)$. Define the profiles
$h,k:\R\to[0,\infty)$ to vanish on $(-\infty,0]$, to be continuous everywhere, and to have the following
constant slopes on the indicated interiors,
\[
\begin{array}{c|cc}
 \text{interval}&h'&k'\\ \hline
 P_0 & b_1 & b_1\\
 I_j,\ j\text{ odd} & b_{j+1} & b_j\\
 I_j,\ j\text{ even} & b_j & b_{j+1}\\
 P_j,\ 1\le j\le m & b_{j+1} & b_{j+1}.
\end{array}
\]
Each profile enters $I_j$ with slope $b_j$ and leaves it with slope $b_{j+1}$, so both are convex. On odd
intervals $h-k$ increases by $\delta_0$, on even intervals it decreases by $\delta_0$, and on connectors
it is constant. Thus $0\le h-k\le\delta_0$ everywhere, and on $P_j$ the difference equals $\delta_0$ for
odd $j$ and zero for even $j$. Across each odd interval $h$ increases by $b_{j+1}\ell_j\ge\delta_0$, so
$h(r_m)\ge(m+1)\delta_0/2$.

We now use the horizontal coordinate $x$ to choose which profile is active at height $y$. At $x=0$ the
first profile dominates, and a horizontal displacement of at least $\delta_0$ makes the second profile,
shifted by $|x|$, active. To also enforce a minimum at height zero, for $(x,y)\in\R^2$ define
\begin{equation}
 f(x,y)\coloneqq\max\{h(y),\, |x|+k(y),\, -\beta y\}.
 \label{gd:lb:objective}
\end{equation}
Each convex profile is the maximum of the affine extensions of its pieces, and $|x|=\max\{x,-x\}$, so $f$
is a finite maximum of affine functions. Their slopes have Euclidean norm at most $\sqrt{1+\beta^2}<2$.
The bounds $f(x,y)\ge|x|$, $f(x,y)\ge-\beta y$, and $f(x,y)\ge h(y)$, together with $h(y)\to\infty$ as
$y\to\infty$, give coercivity. For $y>0$ the first profile is positive, for $y<0$ the term $-\beta y$ is
positive, and $f(x,0)=|x|$, so the origin is the unique minimizer.

Let $g_x,g_y$ denote the components of the selector $g$. The relevant horizontal coordinates form the set
$\mathcal H\coloneqq\{0,\delta_0,\delta_0-1\}$. Call $0$ the inner coordinate and the other two the outer
coordinates. For $y>0$ let $s_h(y)$ and $s_k(y)$ be the interior slopes of the profiles on a long
interval, and on each closed connector use its common slope, including its endpoints. Both choices are
valid subgradients of the respective profiles. Let $\operatorname{sgn}$ denote the sign function with
$\operatorname{sgn}0=0$. For positive heights the default choices are $g(0,y)\coloneqq(0,s_h(y))$ at the inner coordinate and
$g(x,y)\coloneqq(\operatorname{sgn}x,s_k(y))$ at an outer coordinate. These choices keep the horizontal
coordinate in the same class. We change the choice at two kinds of ties on the connectors so that the
next interval is entered in the class with the flatter profile.
\begin{itemize}
 \item On a connector with $h=k$ and common slope $s$, choose $g(0,y)\coloneqq(-\delta_0,s)$. The term $|x|+k(y)$ is active there and $-\delta_0\in[-1,1]=\partial|\cdot|(0)$.
 \item On a connector with $h-k=\delta_0$ and common slope $s$, choose $g(\delta_0,y)\coloneqq(\delta_0,s)$. Both profiles are active there and the vector is the convex combination $(1-\delta_0)(0,s)+\delta_0(1,s)$.
\end{itemize}
Since $h\ge k$ at the inner coordinate and $|x|\ge\delta_0\ge h-k$ at either outer coordinate, the
ordinary choices also select active branches. At height zero choose $g(0,0)\coloneqq0$ and
$g(x,0)\coloneqq(\operatorname{sgn}x,0)$ at an outer coordinate. At a negative height on one of these
three lines, choose $(0,-\beta)$ when $-\beta y>|x|$ and $(\operatorname{sgn}x,0)$ otherwise. These
choices follow from $h=k=0$ there. Everywhere else, fix an order of the finitely many affine functions
representing $f$ and select the slope of the first active one. This defines a Borel selection on the
entire plane with $|g_y|\le\beta$.

The horizontal update $x\mapsto x-g_x(x,y)$ preserves $\mathcal H$. At positive heights the ordinary rule
fixes zero and interchanges the outer coordinates. At a positive height on an even connector, one update
from the inner coordinate gives $\delta_0$ and outer coordinates stay outer. On an odd connector, one
update sends $\delta_0$ to zero, two send $\delta_0-1$ to zero, and zero stays fixed. The rules at height
zero and at negative heights also preserve $\mathcal H$.

We now check that a trajectory cannot cross a connector too quickly for this switch to occur. For $1\le
j\le m$, consider a trajectory starting at the origin if $j=1$, or starting with horizontal
coordinate in $\mathcal H$ and height strictly below $r_{j-1}$ if $j\ge2$. Suppose its vertical
increments have absolute value at most $1/2$. Until its first visit to a height at least $r_j$, its
selected vertical subgradient is at most $b_j$, with equality in the interior of $I_j$. To prove this,
let $c$ be the left endpoint of the connector entering $I_j$, so $c=0$ for $j=1$ and $c=r_{j-1}$
otherwise. At any time $t$ with height in $(a_j,r_j)$ before that visit, let $s<t$ be the last time with
height at most $c$. The initial condition guarantees such an $s$. Since the connector has width ten,
times $s+1$ and $s+2$ are strictly before $t$ and both heights belong to $(c,a_j)$. These two updates
place the horizontal coordinate in the outer class for odd $j$ and in the inner class for even $j$. The
entering connector and the interior of $I_j$ preserve that class until $t$, even on repeated visits, so
the selected vertical slope there is $b_j$. At earlier positive heights all selectable slopes are at most
$b_j$, and at nonpositive heights the selected vertical component is zero or $-\beta$, proving the
assertion.

\subsection{Two estimates for walks with sign increments}
\label{gd:lb:random-walks}

The vertical walk in each stage must do two things. It must reach a target despite a small downward
drift and, after its first visit, it must be back at the target at the prescribed final time. We derive
a separate estimate for each event.

For the hitting estimate, let $(S_t)_{t\ge0}$ be a walk from zero with independent equiprobable increments in $\{-1,1\}$. For
integers $N,A\ge1$, reflecting a path after its first visit to $A$ pairs paths that hit $A$ and end below
$A$ with paths ending above $A$. The interval $[-A,A)$ contains $A$ integers of the parity of $N$. Let
$c_N\coloneqq2^{-N}\binom{N}{\lfloor N/2\rfloor}$ be the largest point mass after $N$ steps, with
$c_0\coloneqq1$. For integers $k\ge0$, $c_{2k}=\prod_{i=1}^k(1-1/(2i))\le(k+1)^{-1/2}$ by the
harmonic-sum bound, and $c_{2k+1}=(2k+1)c_{2k}/(2k+2)$. Thus $c_N\le\sqrt{2/(N+1)}$, and the reflection
gives
\[
 \Pp\left[\max_{0\le t\le N} S_t<A\right]=\Pp[-A\le S_N<A]\le A c_N\le4A/\sqrt N.
\]
For a symmetric walk from one, write $\Pp_1^{\mathrm{sym}}$ for its law and $\tau_0$ for its first visit to
zero. By symmetry, for an integer $u\ge0$, the probability $\Pp_1^{\mathrm{sym}}[\tau_0>u]$ equals the
probability of staying below one from zero for $u$ steps. The same reflection with $A=1$ bounds this by
$c_u\le2/\sqrt{u+1}$, including $u=0$.

We next allow the small downward drift that comes from the selected subgradient. For $0\le\delta\le1/2$,
let $\Pp_\delta$ denote the law of $N$ independent signs of mean $-\delta$, so
$\Pp_0$ is the symmetric law. Under $\Pp_0$, denote the signs by $s_1,\ldots,s_N$ and expectation by
$\E_0$. The likelihood ratio $D\coloneqq\prod_{i=1}^N(1-\delta s_i)$ has $\E_0[D]=1$ and
$\E_0[D^2]=(1+\delta^2)^N$. If $N\delta^2\le1$, then $e^{N\delta^2}-1\le4N\delta^2$. For any event $H$
determined by the signs, Cauchy-Schwarz therefore gives
\[
 |\Pp_\delta[H]-\Pp_0[H]|
 \le\tfrac12\E_0[|D-1|]\le\tfrac12\sqrt{e^{N\delta^2}-1}\le\delta\sqrt N .
\]

After the first hit we measure distance below the target. The next chain bounds this distance from
above. Let $0<\delta\le1/2$ and let $(K_t)_{t\ge0}$ be a chain on the nonnegative
integers, starting at zero. From zero it moves to one with probability $\delta$ and otherwise stays, and
from a positive state it moves up or down by one with respective probabilities $p_+\coloneqq(1+\delta)/2$
and $p_-\coloneqq(1-\delta)/2$. For an integer $N\ge1$ we claim that
\begin{equation}
 \Pp[K_N>0]\le2\delta^2N+4\delta\sqrt N.
 \label{eq:away}
\end{equation}
To be away at time $N$, the chain must have left zero and its last such excursion must still be in
progress. We first bound the duration of an excursion starting from one, accounting also for the chance
that it never returns. Let $r\coloneqq p_-/p_+<1$. For integers $1\le k<B$, the function $k\mapsto(r^k-r^B)/(1-r^B)$ solves the
two-boundary difference equation with value one at zero and zero at $B$. It therefore gives the chance of
hitting zero before $B$ from $k$. Letting $B\to\infty$ shows that the probability of ever hitting zero
from $k$ is $r^k$. Conditional on hitting zero, the probability of an upward move from $k$ is
$p_+r^{k+1}/r^k=p_-$. Coupling this conditioned walk with a symmetric walk from the same point using a
common uniform at each update shows that the conditioned walk hits zero no later than the symmetric walk.
Thus, writing $\Pp_1$ for the law starting from one and $\tau_0$ for the first visit to zero, for
integers $u\ge0$ we have $\Pp_1[\tau_0>u]\le1-r+\Pp_1^{\mathrm{sym}}[\tau_0>u]\le2\delta+2/\sqrt{u+1}$. On
$\{K_N>0\}$ there is a last departure from zero. For $0\le s<N$, the probability that this departure
occurs at time $s$ is $\Pp[K_s=0]\delta\Pp_1[\tau_0>N-s-1]$. Summing and using $\sum_{i=1}^N
i^{-1/2}\le2\sqrt N$ proves \eqref{eq:away}. In particular, for $0<b\le1/4$ with $b\sqrt N\le1$, taking
$\delta=2b$ bounds the probability by $16b\sqrt N$.

\subsection{The noise law and the probability of crossing each interval}
\label{gd:lb:stages}

For the unit-step trajectory denote the coordinates by $(U_t,Y_t)\coloneqq Z_t$, and let $a\coloneqq1/2$
be the vertical lattice spacing. For $1\le j\le m$, let $q_j\coloneqq r_j-4$ be the target height of
stage $j$ and let $N_j\coloneqq Q\ell_j^2$ be its duration. Let $s_0\coloneqq0$,
$s_j\coloneqq\sum_{i=1}^j N_i$, and $H(M,m)\coloneqq s_m$. Stage $j$ runs from $s_{j-1}$ to $s_j$ and
succeeds if $Y_{s_j}=q_j$. Its starting height is zero for $j=1$, and $q_{j-1}$ for $j>1$ if the
preceding stages succeeded. Each target belongs to $a\mathbb Z$ and lies four units inside the right endpoint of its long interval.
This lets us use the smaller slope at the target itself, and stopping four units short costs at most
$4\beta$ in profile value.

We prescribe the height increment directly. Below the target we want a nearly symmetric walk. At the target we allow only a hold or one downward
step. Centering then fixes the two probabilities. Take independent uniform variables on $[0,1]$, one for
each update, and let $\cF_t$ be the $\sigma$-field
generated by the first $t$ variables. As long as every completed stage has succeeded, at each $s_{j-1}\le
t<s_j$ let $v\coloneqq g_y(Z_t)$ and choose a vertical increment $D_t$ using the next uniform. If
$Y_t<q_j$, let $D_t$ equal $a$ or $-a$ with respective conditional probabilities $(1-v/a)/2$ and
$(1+v/a)/2$. If $Y_t=q_j$, let $D_t$ equal $0$ or $-a$ with respective conditional probabilities
$1-b_j/a$ and $b_j/a$. The lattice prevents an overshoot before hitting $q_j$, and the rule at $q_j$
prevents one afterwards. Thus the trajectory remains below $r_j$ throughout an active stage. By the
connector argument in Appendix~\ref{gd:lb:geometry}, the selected slope at $q_j$ is $b_j$. Since
$|v|\le\beta<a$, these are valid probabilities and $\E[D_t\mid\cF_t]=-v$ in both cases. The vertical
innovation $\xi_{t+1}\coloneqq(0,D_t+v)$ is therefore centered and has norm at most $a+\beta<1$. After
the first failed stage use zero innovations. This defines the process on every path.

\begin{lemma}[Success of one planar stage]\label{gd:lb:stage-success}
For each $1\le j\le m$, on the event that all preceding stages succeeded, $\Pp[Y_{s_j}\ne
q_j\mid\cF_{s_{j-1}}]\le1/(8m)$.
\end{lemma}

\begin{proof}
Fix a stage $j$ and condition on its starting history, with all earlier stages successful. We first
bound the probability of missing the target. We then bound the probability of finishing below it after a
visit. The first target is at height $w+\ell_1-4<2\ell_1$, and each later target is $\ell_j+w\le2\ell_j$ above the
preceding one. The initial distance in units of $a$ is therefore an integer $A\le4\ell_j$. Until the
target is reached the connector argument gives $v\le b_j$, so common uniforms bound the increments in
lattice units from below by independent signs with upward probability $(1-2b_j)/2$. In particular, if
this comparison walk reaches level $A$ within $N_j$ steps, the actual trajectory has reached the target.
The stage duration was chosen so that the initial distance and the accumulated drift, each measured in
units of the typical fluctuation $\sqrt{N_j}$, are small. Indeed,
\[
 \frac{\ell_j}{\sqrt{N_j}}=\frac1{\sqrt Q},\qquad
 b_j\sqrt{N_j}=\frac{\delta_0\sqrt Q}{Q-1}\le\frac{2\delta_0}{\sqrt Q}
 =\frac1{1024m}.
\]
Thus $N_j(2b_j)^2\le1$, and the symmetric hitting estimate and biased-sign comparison in
Appendix~\ref{gd:lb:random-walks} bound the probability of failing to reach $q_j$ by
$4A/\sqrt{N_j}+2b_j\sqrt{N_j}\le1/(16m)+1/(512m)$.

For the second requirement we can allow arbitrarily many departures and returns, since success depends
only on the height at $s_j$. Suppose the first visit to $q_j$ occurs at time $\tau\le s_j$. Conditional on $\cF_\tau$, let $R\coloneqq
s_j-\tau$ be the remaining number of updates. For $0\le i\le R$, define the distance
$D_i'\coloneqq(q_j-Y_{\tau+i})/a$, a nonnegative integer. From zero it moves to one with probability
$\delta\coloneqq2b_j$. From a positive distance at time $\tau+i$, the slope bound $g_y(Z_{\tau+i})\le b_j$ gives probability at
most $(1+\delta)/2$ of moving one step farther down. For the chain in \eqref{eq:away}, the transition laws from positive states
are translates of one another. From zero and from one the possible next states are respectively $\{0,1\}$
and $\{0,2\}$, with $\delta\le(1+\delta)/2$. Its kernel is therefore increasing in stochastic order.
Starting both processes at zero and using the same uniforms at each transition keeps the chain above the
actual distance, conditional on the full past. If $R>0$, \eqref{eq:away} gives
$\Pp[D_R'>0\mid\cF_\tau]\le16b_j\sqrt R\le1/(64m)$, and for $R=0$ the distance is zero. Adding the
chances of missing the target and of finishing below it after a visit gives conditional failure
probability at most $41/(512m)<1/(8m)$.
\end{proof}

\subsection{The terminal error at a prescribed horizon}
\label{gd:lb:finish-plane}

\begin{proof}[Proof of Lemma~\ref{gd:thm:bounded-plane}]
At unit stepsize, Lemma~\ref{gd:lb:stage-success} bounds the probability of a first failure at each
stage, conditional on the success of its predecessors, by $1/(8m)$. Summing the probabilities of a first failure gives success of all stages
with probability at least $7/8$ at time $H(M,m)$. On this event the terminal height is $r_m-4$, so
\begin{equation}
 h(r_m-4)\ge h(r_m)-4\beta\ge\delta_0\bigl((m+1)/2-1/3\bigr)\ge m/16.
 \label{gd:lb:height-gain}
\end{equation}
This estimate is uniform in $M\ge1000$ and bounds the profile itself, which we will use when combining
several planar constructions.

The gain is now proportional to the number of stages. To express it in terms of the prescribed horizon,
we estimate the time required for $m$ stages and take the largest admissible $m$. For the planar
assertion take $M=1000$ and abbreviate $H_m\coloneqq H(1000,m)$. The durations give
$H_m=M^2 Q\sum_{i=0}^{m-1}Q^{2i}$, so $M^2 Q^{2m-1}\le H_m\le mM^2 Q^{2m-1}$. As $m$ ranges over the
positive odd integers, $H_m$ increases strictly from $H_1=n_0-1$ to infinity. Using $(2m-1)\log Q\ge
m\log65536$ for the lower bound and $\log M^2<14$ and $\log Q<12+2\log m$ for the upper bound yields
$11m<\log H_m\le40m\log(em)$.

For any integer $n\ge n_0$, choose the largest odd $m$ with $H_m\le n-1$. Since $g(0)=0$, zero innovations
keep the trajectory at zero for $n-1-H_m$ updates, after which the stages end exactly at time $n-1$. We have
$n-1<H_{m+2}$, $m+2\le3m$, and $\log\log(n-1)\ge\log\log H_m>\log(em)\ge1$. Hence $\log(n-1)\le40(m+2)\log(e(m+2))\le360m\log\log(n-1)$. Since $f(x,y)\ge h(y)$, on the all-success event
\eqref{gd:lb:height-gain} gives unit-step terminal error at least $m/16\ge A_n/5760$. This event has
probability at least $7/8$.

For any $\rho>0$, scale each unit-step point by $\rho$, use objective $z\mapsto\rho f(z/\rho)$ and
selector $z\mapsto g(z/\rho)$, and keep the same innovations. The subgradient inequality gives the stated
selector for the scaled objective. Its Lipschitz bound, minimizer, coercivity and finite affine
representation persist, and the scaled recursion has stepsize $\rho$ and terminal value $\rho f(Z_{n-1})$ in
terms of the original trajectory. This multiplies the threshold in the probability bound by $\rho$. Since
$(7/8)/5760\ge1/10000$ and the objective is nonnegative, it also gives the mean bound in
\eqref{gd:lb:plane-prob}.
\end{proof}

For the normalized class at stepsize $\eta$, take $\rho=\eta/2$ and divide the objective and selector by
two. Set the noise at update $t+1$ in \eqref{eq:ssgm} to $-\xi_{t+1}/2$. This realizes the same
trajectory, with Lipschitz constant less than one and noise norm at most $1/2$. Starting at the unique
planar minimizer, it gives
\begin{equation}
 \Pp[\Delta(X_n)\ge\eta A_n/23040]\ge7/8,\qquad
 \E[\Delta(X_n)]\ge\eta A_n/40000 .
 \label{gd:eq:normalized-plane}
\end{equation}
Appending zero coordinates gives the same lower bounds in every dimension $d\ge2$. The next two constructions retain the planar order under additional restrictions on the oracle or the domain.

\subsection{A single state-dependent noise law in the plane}
\label{app:homogeneous}

In the construction above, a point visited in two different stages can receive different noise laws.
We now realize the same lower order with one law at each point. We do this on the full plane, where
we can perturb the finitely many states of the original construction. The law and objective may
still depend on the prescribed $n$ and $\eta$.

Fix $n\ge n_0$ and $\eta>0$, and take the normalized planar instance of
\eqref{gd:eq:normalized-plane}. Its possible histories form a finite tree because each update has
finitely many outcomes. For a node $v$ let $x_v$ be its state and let $g_v\in\partial f(x_v)$ be its
chosen subgradient, choosing any such vector at a leaf. At the root $o$ we have $x_o=g_o=0$.
For each child $i$ of a nonterminal node $v$, let $p_{vi}$ be its transition probability and let
$w_{vi}$ be its noise. Writing $\mathcal C(v)$ for the children, we have
$x_i=x_v-\eta(g_v+w_{vi})$, $\sum_{i\in\mathcal C(v)}p_{vi}w_{vi}=0$, and $\|w_{vi}\|\le1/2$.

Let $V\coloneqq\eta A_n/40000$ be the lower bound on the original expected error. We use the
Moreau envelope $M_\lambda$ defined in Appendix~\ref{app:tracking}, now for the planar objective
with minimum zero. For every $x\in\R^2$ and $g\in\partial f(x)$, proximal optimality gives
$P_\lambda(x+\lambda g)=x$. The envelope has gradient
$\nabla M_\lambda(z)=(z-P_\lambda z)/\lambda$, so at $x+\lambda g$ its gradient is exactly $g$
and its value is $f(x)+\lambda\|g\|^2/2$. The envelope is convex and $1$-Lipschitz with
$1/\lambda$-Lipschitz gradient. It has the same minimizers as $f$, since the two nonnegative
terms in its definition vanish together. It is also coercive, since bounded envelope values bound
both the value of $f$ at the proximal point and the distance to that point.

We use these identities to assign distinct points to distinct histories. Let
$\alpha\coloneqq1/(8(n-1))$ and choose $\lambda>0$ small enough that $\lambda^2\le V/8$ and
$\lambda+2\lambda/\eta+2\lambda^2/\eta\le\alpha/2$. At each node let
$y_v\coloneqq x_v+\lambda g_v+e_v$, where $e_o=0$ and $\|e_v\|\le\lambda^2$. Since each ball of
possible perturbations contains infinitely many points and the tree is finite, we can choose the
$y_v$ to be pairwise distinct. The Lipschitz bounds for the envelope and its gradient then give
\[
 M_\lambda(y_v)\ge f(x_v)-\lambda^2,\qquad
 \|\nabla M_\lambda(y_v)-g_v\|\le\lambda .
\]

To go from $y_v$ to $y_i$ in one step, the required noise is
$\widetilde w_{vi}\coloneqq(y_v-y_i)/\eta-\nabla M_\lambda(y_v)$. By the original update,
$\|\widetilde w_{vi}-w_{vi}\|\le2\lambda/\eta+2\lambda^2/\eta+\lambda\le\alpha/2$.
Thus $m_v\coloneqq\sum_{i\in\mathcal C(v)}p_{vi}\widetilde w_{vi}$ has norm at most $\alpha/2$.
At $y_v$, choose the original child with probability $1-\alpha$, distributed according to its
probabilities $p_{vi}$, and use noise $\widetilde w_{vi}$. With the remaining probability $\alpha$,
use noise $-(1-\alpha)m_v/\alpha$. This law is centered and supported on the ball of radius
$1/2+\alpha/2<1$. At all other points use zero noise. Because the node labels are distinct, this
defines a single Borel noise kernel on the plane.

Use fresh independent Bernoulli variables for the choice between the original branches and the
correction. With probability $(1-\alpha)^{n-1}\ge7/8$ no correction occurs, and conditional on
this event the labels follow exactly the original tree law. Since the terminal value is nonnegative
also on the other paths, for the resulting last iterate $Y_n$ we obtain
\[
 \E[M_\lambda(Y_n)]\ge(1-\alpha)^{n-1}(V-\lambda^2)
 \ge(7/8)^2 V\ge\frac{\eta A_n}{80000}.
\]
Thus a single state-dependent law realizes the planar order from the unique minimizer of a
differentiable objective. The construction allows its gradient Lipschitz constant to depend on $n$
and $\eta$.

\subsection{Keeping the trajectory in a fixed ball}
\label{gd:lb:localization}

We first show that the event of a large final error still has positive probability when we also require
the whole path to stay in a fixed ball. On this event projection changes no update. We will then handle
the paths that leave the ball to obtain an ordinary, unprojected trajectory as well. For $R>0$ let
$B_R\coloneqq\{z\in\R^2:\|z\|\le R\}$. Fix $n\ge n_0$ and $0<\rho\le 1/\sqrt{n-1}$ and use the
planar instance of Lemma~\ref{gd:thm:bounded-plane}. Let $\tau$ be the first $t\in\{0,\ldots,n-1\}$ for
which $\|Z_t\|>4$, taking $\tau=\infty$ if no such $t$ exists. Since the origin minimizes $f$, the
subgradient inequality gives $\langle z,g(z)\rangle\ge f(z)-f(0)\ge0$. Expanding the squared norm,
conditioning on $\cF_t$, and using $\|g(Z_t)\|<2$ and $\|\xi_{t+1}\|\le1$ gives
$\E[\|Z_{t+1}\|^2\mid\cF_t]\le\|Z_t\|^2+5\rho^2$. Multiplying the increment bound by $\1_{\{t<\tau\}}$
and summing yields $\E[\|Z_{\min\{n-1,\tau\}}\|^2]\le5(n-1)\rho^2$. On exit the squared norm exceeds $16$, so
$\Pp[\tau<n]\le5/16$. Together with \eqref{gd:lb:plane-prob} this gives
\[
 \Pp[f(Z_{n-1})\ge\rho A_n/5760,\ \tau\ge n]\ge7/8-5/16=9/16 .
\]

On $B_4$, run the projected updates with the same selector and innovations. The innovations remain
conditionally centered under the original filtration, and on $\{\tau\ge n\}$ every proposed update lies in
$B_4$, so the projected and planar trajectories coincide. On $B_8$, run ordinary updates with innovations
$\1_{\{t<\tau\}}\xi_{t+1}$. Since the indicator is $\cF_t$-measurable, centering and the noise bound
hold. This trajectory agrees with the planar one up to and including exit, and through time $n-1$ on
$\{\tau\ge n\}$. At exit its norm is at most $4+3\rho$. Thereafter the noise vanishes, and the subgradient
inequality bounds the increase of its squared norm by $4\rho^2$ per step. Thus all post-exit squared
norms are at most $(4+3\rho)^2+4(n-1)\rho^2\le53<8^2$. This gives an ordinary trajectory strictly inside
$B_8$ on every path. In either construction the terminal error is at least $\rho A_n/5760$ with
probability at least $9/16$, hence has expectation at least $\rho A_n/10240$. In particular
$\rho=1/\sqrt{n-1}$ gives order $(\log n)/(\sqrt n\log\log n)$ on fixed balls. Since the error on a domain of
finite diameter is bounded by $L$ times that diameter, the fixed-ball estimate cannot grow indefinitely
at a fixed positive stepsize.

It remains to handle every prescribed stepsize on the same ball. We scale the objective and noise
together so that the actual displacement stays within the range of the exit estimate, then use separate
examples when that scaling makes the error too small. Let $\alpha\coloneqq\tfrac12\min\{1,1/(\eta\sqrt{n-1})\}$ and $\rho\coloneqq\eta\alpha$.
Since $\rho\le1/(2\sqrt{n-1})$, we can use the preceding construction at this stepsize. For stepsize $\eta$,
scale the objective and selector by $\alpha$ and the noise by $-\alpha$ when adopting the sign convention
of \eqref{eq:ssgm}. The trajectory on either ball is unchanged and its expected error is at least
$\alpha\rho A_n/10240\ge(A_n/40960)\min\{\eta,1/(\eta n)\}$.

Two further instances on the same domains cover small and large steps. Write $e_1,e_2$ for the coordinate
unit vectors, and $\|x\|_\infty$ for the maximum absolute coordinate. For small steps let
$a\coloneqq\min\{1/2,1/(2\eta n)\}$ and take $f(x)\coloneqq a\|x\|_\infty$ with zero noise, starting at
$e_1$. Selecting $ae_1$ along the positive first coordinate axis gives, for $1\le k\le n$,
$X_k=(1-(k-1)\eta a)e_1$. These iterates have norm between $1/2$ and one, which verifies the selection
and gives terminal error at least $\tfrac14\min\{1,1/(\eta n)\}$. For large steps take
$f(x)\coloneqq\|x\|_\infty$, start at zero, select zero there, and choose zero noise until the last
update. At this update use a uniform sign times $\min\{1,1/\eta\}e_2$. The last error is
$\min\{\eta,1\}$, and all points remain in $B_1$.

Each instance has $D_1\le1$, $L\le1$ and noise norm at most one. Let $c_0\coloneqq1/40960$. If $\eta\le
n^{-1/2}$ and $\eta n\ge1$, the geometric mean of the planar and small-step bounds is
$(\sqrt{c_0}/2)\sqrt{A_n/n}$. If $\eta n<1$, then $A_n\le n$ and the small-step bound $1/4$ already
suffices. For $n^{-1/2}\le\eta\le1$, the geometric mean of the planar and large-step bounds is
$\sqrt{c_0}\sqrt{A_n/n}$, and for $\eta\ge1$ the large-step bound is one. Thus in either domain setting,
for each $\eta$ there is an instance with expected error at least $(1/(128\sqrt{10}))\sqrt{A_n/n}$.

For $d>2$, extend planar objectives to depend only on their first two coordinates and append zeros to the
subgradients and noises. Projection onto a Euclidean ball preserves this plane. Finally, for
$R\in\{4,8\}$ let $s\coloneqq1/R$. To work on the unit ball with stepsize $\eta$, construct the instance
on $B_R$ at stepsize $\eta/s$, scale its points by $s$, and use the objective $x\mapsto s f(x/s)$.
Subgradient and noise amplitudes are unchanged, projection commutes with the scaling, and the initial
distance is at most $s$. The error is multiplied by $s\ge1/8$. Since $A_n$ has order
$(\log n)/\log\log n$ as $n\to\infty$, for $D=1$ and every $d\ge2$ this gives an optimized lower
bound of order $\sqrt{\log n/(n\log\log n)}$ on the unit ball, both for projected updates and for
ordinary iterates that stay strictly inside the ball. For each fixed $d\ge2$, the uniform upper bound in Theorem~\ref{thm:main}
at $p=\infty$ gives the same optimized order for these bounded-noise subclasses.

\section{A lower bound depending jointly on dimension and horizon}
\label{gd:app:joint}

We prove that there is a universal $c_\star>0$ such that for every $d\ge2$, $n\ge n_0$ and $\eta>0$,
there is a normalized instance on $\R^d$ starting at a minimizer, with bounded centered noise, for which
\begin{equation}
 \E[\Delta(X_n)]\ge c_\star\eta q_{n,d}.
 \label{gd:eq:joint-lower}
\end{equation}
Its objective will be a finite maximum of affine functions and its selected subgradient a fixed Borel
function of the point. For $d=2$, $\log n\asymp\log(n-1)$ and $\log(e+(\log n)/2)\asymp\log\log(n-1)$ throughout this range of horizons,
so \eqref{gd:eq:normalized-plane} suffices. For $d\ge3$, we run a planar construction on one coordinate pair at a time. Two requirements govern the
construction. The terminal objective should count the gains from all pairs, yet the Lipschitz bound must
be independent of how many pairs we use. Moreover, earlier pairs must retain their gains while later
pairs run. We first build an objective and a schedule meeting these requirements, then use one extra
coordinate to make the active pair identifiable from the current point.

\subsection{A planar construction for one coordinate pair}

Let $m\ge1$ be odd, let $M\ge1000$ be an integer, and use the profiles, slopes and intervals from
Appendix~\ref{gd:lb:geometry}. We first adapt one planar construction for use inside the common objective. At negative heights we give
the two profiles the same slope, which will allow an inactive pair to contribute a function of height
alone. For $y\in\R$ define $\phi(y)\coloneqq h(y)+\beta(-y)_+$ and
$\kappa(y)\coloneqq k(y)+\beta(-y)_+$. Both are convex and $\beta$-Lipschitz, satisfy
$0\le\phi-\kappa\le\delta_0$, and have slope $-\beta$ on the negative half-line. Use the planar objective
$F_M(x,y)\coloneqq\max\{\phi(y),|x|+\kappa(y)\}$. At positive heights it agrees with
\eqref{gd:lb:objective}, so keep the selector already constructed there. For $y<0$ use
$(\operatorname{sgn}x,-\beta)$, and for $y=0$ use $(\operatorname{sgn}x,0)$. At negative heights
$F_M(x,y)=|x|-\beta y$, and at zero height $|x|+\kappa(y)$ is active with $0\in\partial\kappa(0)$, so
both choices give subgradients.

Let $K_m\coloneqq Q\sum_{j=0}^{m-1}Q^{2j}$. At unit stepsize the $m$ stages from the origin take
$H(M,m)=M^2K_m$ updates. Since at nonpositive heights the horizontal update fixes zero and interchanges
the outer coordinates, the horizontal coordinate stays in $\mathcal H$ throughout, including after a
failed stage when innovations vanish. At nonpositive heights the selected vertical slope is at most zero. The connector argument therefore
still bounds the slope by $b_j$ before the target and gives equality at the target, as required in the
proof of Lemma~\ref{gd:lb:stage-success}. The two quantities in that probability estimate, $\ell_j/\sqrt{N_j}=1/\sqrt Q$ and
$b_j\sqrt{N_j}=\delta_0\sqrt Q/(Q-1)$, are independent of $M$. Thus for every $M\ge1000$, with
probability at least $7/8$ the terminal vertical coordinate is $r_m-4$, and \eqref{gd:lb:height-gain}
gives $h(r_m-4)\ge m/16$.

\subsection{A common objective and a schedule of updates}

Let $q$ be the number of pairs, an integer between one and $\lfloor(d-1)/2\rfloor$. We will give a finished pair no further vertical noise. Its vertical coordinate can then move by at most
its slope bound per update, and the profile value can fall by at most the square of that bound. Since the slope bound decreases as $1/M$, we choose
the length parameter for each pair in terms of how many updates will follow it. This explains why the
schedule is defined backward. Let $S_{q+1}\coloneqq0$. For
$i=q,q-1,\ldots,1$ define
\[
 M_i\coloneqq\left\lceil1000\sqrt{1+S_{i+1}}\right\rceil,\qquad
 H_i\coloneqq K_m M_i^2,\qquad S_i\coloneqq H_i+2+S_{i+1}.
\]
Here $H_i$ counts the planar updates of pair $i$, the next two return its horizontal coordinate to zero,
and $S_i$ counts all updates from the start of pair $i$ onward. Denote the corresponding profiles and
slope bound by $\phi_i,\kappa_i,h_i$ and $\beta_i$. To count all vertical gains, we include their sum as one branch of the objective. In each additional
branch we replace exactly one vertical profile by the term that activates its horizontal coordinate.
Thus only one horizontal slope contributes to the Lipschitz bound of a branch. For horizontal and
vertical coordinate vectors $x,y\in\R^q$, let
\[
 F(x,y)\coloneqq\max\left\{
 \sum_{j=1}^q\phi_j(y_j),\
 \max_{1\le i\le q}\left(|x_i|+\kappa_i(y_i)+
       \sum_{\substack{1\le j\le q\\j\ne i}}\phi_j(y_j)\right)\right\}.
\]
The summands are nonnegative and the origin has value zero. All profiles are convex and piecewise affine
with finitely many pieces, so $F$ is a finite maximum of affine functions. In any such affine piece at
most one horizontal slope is nonzero, of magnitude at most one, and the vertical slopes have magnitudes
at most $\beta_1,\ldots,\beta_q$. Since $\beta_i\le2\delta_0/M_i$ and $1+S_i\ge2(1+S_{i+1})$,
\[
 \sum_{i=1}^q\beta_i^2\le(2\delta_0/1000)^2\sum_{a=0}^{q-1}2^{-a}<1 .
\]
It follows that the squared Lipschitz constant of $F$ is at most $1+\sum_{i=1}^q\beta_i^2$.

We now check that the common objective reproduces the planar dynamics during the updates assigned to one
pair. Fix $i\in\{1,\ldots,q\}$ and suppose that all horizontal coordinates except possibly $x_i$ vanish.
Since $\kappa_j\le\phi_j$, every indexed branch with $j\ne i$ is bounded by the first branch at this
point. The function $F_{M_i}(x_i,y_i)+\sum_{1\le j\le q,\,j\ne i}\phi_j(y_j)$ is a global lower bound on
$F$ and equals it here. Any subgradient of this lower bound at the point is consequently a subgradient of
$F$. In pair $i$ use the planar selector. In every other pair use zero horizontal component and any fixed
one-sided slope of $\phi_j$ in the vertical component, choosing zero at $y_j=0$. In particular, pairs
that have yet to start stay at the origin.

Before switching to the next pair, we must return the current horizontal coordinate to zero on every
path. A single centered innovation cannot in general cancel a known nonzero displacement. The first of
two updates rounds the tentative coordinate to a value in $\{-1,0,1\}$ with the same conditional mean,
and the second uses a valid subgradient to send that value to zero. After the $H_i$ planar updates of
pair $i$, let $v$ denote the tentative horizontal coordinate from one
more update with the planar selector and no noise. We have $v\in\mathcal H$ and $|v|\le1$. Given the
history, choose a variable $B$ equal to $\operatorname{sgn}v$ with probability $|v|$ and to zero
otherwise. Use $B-v$ as the only nonzero component of the innovation, in coordinate $x_i$. It is
conditionally centered, has magnitude at most one and leaves $x_i=B\in\{-1,0,1\}$. In the next update use
zero innovation. If $x_i=\pm1$, the indexed branch exceeds the first by at least $1-\delta_0>0$, so its
horizontal slope $x_i$ sends the coordinate to zero. If $x_i=0$, the first branch is active and we choose
zero horizontal slope. For the vertical components choose slopes of the branch selected in each case.
After these two updates all horizontal coordinates vanish on every path, as required to start the next
pair.

Use a fresh independent uniform at each randomized update, with the filtration generated by draws already
used. Conditional on the full history at the start of pair $i$, that pair starts from the origin and
follows the planar rule for $H_i$ updates. At the end of those updates the profile $h_i$ is nonnegative
and reaches at least $m/16$ with probability at least $7/8$, so its expected value is at least $7m/128$.
In the remaining $2+S_{i+1}$ updates the vertical coordinate receives zero innovation, and its
displacement per step is at most $\beta_i$. Since $h_i$ is $\beta_i$-Lipschitz, the total decrease of
that profile is at most
\[
 (2+S_{i+1})\beta_i^2\le2(2\delta_0/1000)^2<1/1000 .
\]
Since $F(x,y)\ge\sum_{i=1}^q\phi_i(y_i)\ge\sum_{i=1}^q h_i(y_i)$, linearity of expectation gives terminal
expected value at least $q(7m/128-1/1000)\ge qm/20$.

\subsection{Choosing the number of pairs and stages}

The preceding construction gives expected gain proportional to $qm$. We now choose how many pairs and
how many stages per pair fit in the horizon. Reserving one of the $n-1$ updates for the extra coordinate leaves $n-2$ updates for the pairs. Since $M_i\le1001\sqrt{1+S_{i+1}}$,
we have $1+S_i\le(4+1001^2K_m)(1+S_{i+1})$. Also $K_m\le mQ^{2m-1}$, so there is a universal $C_0\ge1$
such that, for all admissible $m,q$, $\log(1+S_1)\le C_0qm\log(em)$. It therefore suffices to allocate $C_0m\log(em)$ per pair on the logarithmic time scale, and each pair
needs two coordinates. We first choose the number of pairs with at least one stage per pair, then use
the remaining time for additional stages. Let $J\coloneqq\log(n-2)$. Once $J\ge4C_0$, choose
\[
 q\coloneqq\min\left\{\left\lfloor\frac{d-1}{2}\right\rfloor,
                      \left\lfloor\frac{J}{4C_0}\right\rfloor\right\}.
\]
Since $d\ge3$, this choice gives $q\ge1$, and we may define $u\coloneqq J/(C_0q)\ge4$. Choose the largest
odd $m$ with $m\log(em)\le u$. By construction $S_1\le n-2$, and since $m\le u$, maximality gives
$u<(m+2)\log(e(m+2))\le3m\log(3eu)$. Thus $qm\gtrsim J/\log(e+J/q)$. If the dimension determines $q$,
then $q\ge d/4$, and if the horizon determines $q$, then $J/q\le8C_0$ once $n$ is sufficiently large. In
both cases $\log(e+J/q)\lesssim\log(e+J/d)$. Since $\log(n-2)\asymp\log n$ and, for $x\ge0$, $\log(e+2x)\le2\log(e+x)$, we obtain
$qm\gtrsim q_{n,d}$ with universal comparison constants.

\subsection{Realizing the schedule by a fixed selector}

So far the prescribed subgradient may depend on which pair is scheduled for the current time. To realize
this schedule with a fixed selector, we make the magnitude of the reserved coordinate decrease by the
same amount on every path. Index the states by the number of completed updates, starting at zero. Let $z$ denote this coordinate, let $\gamma\coloneqq1/(4n)$, and add $\gamma|z|$ to $F$ to
obtain $\widehat F$. Ignore any unused coordinates. The squared Lipschitz constant of $\widehat F$ is at
most $1+\sum_{i=1}^q\beta_i^2+\gamma^2<4$. To start this coordinate from the minimizer with a centered innovation, choose zero subgradient at the
origin and give the first $z$ innovation values $1/2$ and $-1/2$ with equal probability. All paired
coordinates stay at zero on this update. Thereafter give $z$ zero innovation and subgradient $\gamma\operatorname{sgn}z$. At each time $1\le
t<n$ the resulting coordinate $z_t$ has magnitude $1/2-(t-1)\gamma>0$. These distinct magnitudes
determine $t$ from the current point.

Keep all pairs at the origin for $n-2-S_1$ updates, then apply the successive planar and reset updates just
described. On each prescribed level of $|z|$, use the corresponding subgradient rule wherever its
coordinate conditions hold. Only the innovation rule depends on past failures, so the same selector
applies throughout a planar stage. These conditions define Borel sets, the selected vectors are
subgradients there, and at every other point we can choose the slope of the first active affine function
in a fixed ordering. This gives a single Borel selector $g$ for $\widehat F$ with $g(0)=0$. Each
innovation is centered conditionally on the past and has norm at most one, and the nonnegative $z$
summand preserves the lower bound on the terminal value.

Let $V_0,\ldots,V_{n-1}$ be the resulting unit-step trajectory and write $\xi_{t+1}$ for the innovation added
at update $t$. It remains to realize this unit-step process at the prescribed stepsize $\eta$ with a unit Lipschitz
bound. We scale points by $\eta/2$ and both the subgradients and the innovations by $1/2$. Accordingly,
define $f(x)\coloneqq(\eta/4)\widehat F(2x/\eta)$. For
$0\le t<n$ let $X_{t+1}\coloneqq(\eta/2)V_t$, and for $0\le t<n-1$ let $G_{t+1}\coloneqq
g(2X_{t+1}/\eta)/2$ and $W_{t+1}\coloneqq-\xi_{t+1}/2$. Substitution gives \eqref{eq:ssgm} on $\R^d$ with
the filtration at index $t+1$ generated by the first $t$ draws. The objective has Lipschitz constant less
than one, $D_1=0$ and the noise norm is at most $1/2$. Moreover,
\[
 \E[\Delta(X_n)]=(\eta/4)\E[\widehat F(V_{n-1})]\ge\eta qm/80
 \gtrsim\eta q_{n,d}.
\]
This proves \eqref{gd:eq:joint-lower} above a universal threshold $n_1\ge n_0$. For $n_0\le n<n_1$, use
the embedded planar instance in \eqref{gd:eq:normalized-plane}. Since $n_0,n_1$ are numerical, $A_n\ge1$ and $q_{n,d}\le\log n_1$ on
this range, so a universal constant also gives \eqref{gd:eq:joint-lower} at these horizons.

\section{Explicit comparison constants}
\label{app:constants}

The lower and upper comparison constants in Theorem~\ref{thm:main} can be taken to be $1/43004$ and
$2^{22}$, respectively. We verify these values for the factor $b_{n,p,d}$ in \eqref{eq:factor}.
With $a_\eta\coloneqq\min\{D,D^2/(\eta n)\}$ as in Appendix~\ref{app:main}, the prescribed-step bounds are
\begin{equation}\label{eq:numerical-risk}
 \frac{a_\eta+\eta b_{n,p,d}}{43004}
 \le\mathcal R_{n,p}^{(d)}(\eta,D)
 \le 2^{22}(a_\eta+\eta b_{n,p,d}).
\end{equation}

\subsection{Lower bound}

By \eqref{eq:lower-all-steps}, the risk is at least $a_\eta/4$ and at least $\eta\Phi_{n,p}/128$.
It remains to prove that for $d\ge2$ it is also at least $\eta q_{n,d}/43000$, where
$q_{n,d}=(\log n)/\log(e+(\log n)/d)$. Let $u\coloneqq\log(n-1)$. If $u<32768$, then
$q_{n,d}\le\log n<32769$. Take $f(x)=|x_1|$ on $\R^d$, start at the origin with zero subgradient
and use zero noise until the last update. A symmetric unit increment in the first coordinate at
that update gives error $\eta$, as required.

Suppose now that $u\ge32768$. Since $u\le\log n\le u+1$, we have
$q_{n,d}\le(32769/32768)u/\log(e+u/d)$. In the plane, $\log u\ge15\log 2$ gives
$\log(e+u/2)\ge(14/15)\log u$. The planar construction in \eqref{gd:eq:normalized-plane}
therefore gives the desired bound, since its expected error is at least $\eta u/(40000\log u)$ and
$40000(15/14)(32769/32768)<43000$.

For $d\ge3$, we specify the number $q$ of coordinate pairs and the odd number $m$ of stages in
Appendix~\ref{gd:app:joint}. For those choices that fit within the horizon, its construction gives
expected error at least $\eta qm/80$. Let $K_m\coloneqq\sum_{j=0}^{m-1}(256m)^{4j+2}$ be the number
of updates in the $m$ planar stages at unit spatial scale. For $1\le i\le q$, let $S_i$ count the
updates from the start of pair $i$ onward. The schedule gives $1+S_i\le(4+1001^2K_m)(1+S_{i+1})$,
with $S_{q+1}=0$.
Since $1\le K_m\le m(256m)^{4m-2}$, we have
\[
 \log(4+1001^2K_m)\le(32m+4)\log2+(4m-1)\log m\le32m\log(em).
\]
Let $j\coloneqq\log(n-2)$ and choose
$q\coloneqq\min\{\lfloor(d-1)/2\rfloor,\lfloor j/128\rfloor\}$. Let $z\coloneqq j/(32q)\ge4$
and take the largest odd $m$ with $m\log(em)\le z$. Iterating the preceding bound gives
$\log(1+S_1)\le32qm\log(em)\le j$, so the construction fits. To bound its gain, maximality and
$m\le z$ give $z<(m+2)\log(e(m+2))\le3m\log(3ez)$, hence $qm\ge j/(96\log(3ez))$.
If $q=\lfloor(d-1)/2\rfloor$, then $q\ge d/4$, so $z\le u/(8d)$ and
$\log(3ez)\le2\log(e+u/d)$. Since $j\ge u/2$, we obtain $qm\ge u/(384\log(e+u/d))$.
If $q=\lfloor j/128\rfloor$, then $j\ge256$ gives $qm\ge j/256\ge u/512$.
Thus in either case the error is at least $\eta u/(40960\log(e+u/d))$, and
$40960(32769/32768)<43000$ gives the desired lower bound.

The two noise constructions start at a minimizer, and the construction for the dimension term
uses bounded noise. Taking the larger applicable stochastic lower bound therefore gives
$\eta b_{n,p,d}\le43000\,\mathcal R_{n,p}^{(d)}(\eta,D)$ for every $p$ and $d$, including $D=0$.
Adding $a_\eta\le4\mathcal R_{n,p}^{(d)}(\eta,D)$ proves the first inequality in \eqref{eq:numerical-risk}.

\subsection{Upper bound}

For $p\le2$, \eqref{gd:eq:heavy-fixed-step} gives the upper bound $a_\eta+8\eta\Phi_{n,p}$.
For $d=1$, \eqref{eq:normalized-upper} gives $a_\eta+(65/2)\eta\Phi_{n,p}$.
We now keep track of the constants in Appendix~\ref{app:geometric} for $d\ge2$ and $p>2$.
We first use unit steps and subtract $f_\star$ from $f$. As in that appendix, write $A\coloneqq D^2/n$,
$a_N\coloneqq\sqrt A(N/n)^{1/6}$ for $N>0$, and $\varepsilon\coloneqq\min\{p-2,1\}$, with
$\varepsilon=1$ at $p=\infty$. Let $J\coloneqq\log n$ and $w\coloneqq J/d$. If $w\le e^{1024}$,
then $\log(e+w)<1025$ and $H_{n-1}+1/2\le3/2+J\le1027b_{n,p,d}$, where $H_{n-1}$ is the
harmonic number. Thus \eqref{gd:eq:variance-normalized} already suffices.

Suppose $w>e^{1024}$. Let $h\coloneqq\log w$, $r\coloneqq2^{\lfloor(\log_2 w)/8\rfloor}$,
$\rho\coloneqq\log r$, and $P\coloneqq\lfloor\log(n-1)/(2\rho)\rfloor$. For $0\le i\le P$ let
$N_i\coloneqq r^{2(P-i)}$ be the length of suffix $i$. We will use
\begin{equation}\label{eq:numerical-scales}
 J/h\le P\le5q_{n,d},\qquad \rho\le P,\qquad r^4d\log(1+P/d)\le P .
\end{equation}
Indeed, $(h-8)/8\le\rho\le h/8$, $\log(n-1)\ge J/2$, and $h^2\le J$.
Taking the integer part gives $J/h\le P\le4J/(h-8)\le5q_{n,d}$ since $\log(e+w)\le h+1$.
For the last inequality, $r^4\le\sqrt w$, $d/P\le h/w$, and $P/d\le5w$, so the ratio of its
left side to $P$ is at most $2h^2e^{-h/2}\le1$. These choices also give $r>2^{100}$.

\subsubsection{Motion and supporting inequalities on a suffix}

We first bound the two parts of the tracking error $H$ in Appendix~\ref{app:tracking}.
For a window of length at most $k$ starting at time $s$, the sum of squared innovations has
conditional expectation at most $k$, and the martingale in the bound for $H^2$ has conditional
variance at most $k^2$. Doob's inequality therefore gives $\|H/\sqrt k\|_{2,s}\le\sqrt6$.
After truncation the sum of squared innovations has conditional second-moment norm at most
$\sqrt5\,k$. Thus the fourth-moment norm of the normalized tracking error after truncation is at most
$(1+\sqrt5+4)^{1/2}<3$.

For the difference between the original and truncated paths, let $\gamma\coloneqq k^{-\varepsilon/2}$.
At offset $j$ let $d_j$ and $h_j$ be the differences of positions and of selected subgradients,
respectively, as in Appendix~\ref{app:tracking}. Its coupling gives
$\E_s[\|h_j\|^2]\le8j\gamma/k$ and hence
$\E_s[\|d_j\|^2]\le2j\gamma+4j(j-1)\gamma/k$.
The conditional variance of the martingale of inner products in the difference recursion is at most
$2\gamma^2\sum_{j=0}^{k-1}(2j+4j(j-1)/k+8j/k)\le(14/3)k^2\gamma^2$.
The squared changes in noise and subgradients have total conditional expectation at most $6k\gamma$.
Applying Doob's inequality to the recursion bounds the expected maximum squared distance, conditional
on $\cF_s$, by $(6+4\sqrt{14/3})k\gamma<16k\gamma$. Thus in \eqref{eq:tracking-pieces} we may take
$\|U_k\|_{4,s}\le3$ and $\|V_k\|_{2,s}\le4k^{-\varepsilon/4}$.

In the recursion for the proximal displacement $R_N(X_t)$ in Appendix~\ref{app:prox-motion}, the
tracking coefficients sum to at most $12$ and the coefficient of $I_{N,t}$ is at most $39^{2/3}<12$.
Consequently in \eqref{eq:prox-split} we can bound $R_N(X_t)/\sqrt N$ by $U+V+12I_{N,t}$ with
$\|U\|_4\le36$ and $\|V\|_2\le48N^{-\varepsilon/4}$. Using the direct tracking bound
$12\sqrt6<30$ gives $\|R_N(X_t)/\sqrt N\|_2\le30+12I_{N,t}$.

For a suffix of length $N$ starting at $s=n-N$, its maximum displacement $Z_N$ satisfies
$Z_N\le3R_N(X_s)+H$. For integers $N/2\le l\le N$, the definition of $I$ gives
$I_{N,s}\le2a_N$ and $I_{l,s}\le3a_N$, since $s-1+N\ge n/2$ and $s-1+l\ge n/4$.
Thus $\|Z_N/\sqrt N\|_2\le93+72a_N$ and $\|R_l(X_s)/\sqrt l\|_2\le30+36a_N$.
With $B\coloneqq2^{14}$ and $b_N\coloneqq1+a_N^2$, Cauchy-Schwarz bounds both squared norms
by $Bb_N$, since $93^2+72^2<B$. The decomposition of $Z_N/\sqrt N$ has deterministic term
$72a_N$, fourth-moment bound $3\cdot36+3=111$ and second-moment bound
$(3\cdot48+4)N^{-\varepsilon/4}=148N^{-\varepsilon/4}$.

We next apply these estimates to the vectors $v,y$ in Appendix~\ref{app:suffixes}. There
$\|v\|\le2Z_N/\sqrt N$, and $y=m-v$, where the weighted noise sum $m$ has second-moment norm
at most one. After truncation the fourth-moment bound for $m$ is $3$, and the difference has
second-moment norm at most $\sqrt2\,N^{-\varepsilon/4}$. Thus each of $\|v\|$ and $\|y\|$ is
bounded by $U+V+144a_N$ with $\|U\|_4\le256$ and $\|V\|_2\le320N^{-\varepsilon/4}$.
The direct bounds give $\E[\|v\|^2]\le4Bb_N$ and $\E[\|y\|^2]\le9Bb_N$.

We can now evaluate the two supporting errors $C_s,E_s$ in \eqref{eq:geom-support}.
Since $l\ge N/2$ and $\lambda=N/8$, the displacement penalty contributes at most $(4/3)Bb_N$ to
$\E[C_s]$. The squared norms of the proposed updates contribute at most $1$ after division by $2l$,
and the normalized martingale sum has expected absolute value at most $\sqrt{6Bb_N}$. In the expression for $\E[E_s]$,
the displacement, drift, and squared update terms are at most $16Bb_N$, $6Bb_N$, and $4$,
respectively, since the total weight is $\nu\ge N/2$. In the last term, the martingale costs at most
$2\sqrt{Bb_N}$ and summation by parts costs at most $2Bb_N+4$. Hence
$\E[C_s]\le2Bb_N$ and $\E[E_s]\le25Bb_N$.

Over the chosen suffixes, geometric summation gives $\sum_{i=0}^P a_{N_i}^2\le2A$ and
$\sum_{i=0}^P a_{N_i}\le2\sqrt A$. The adjacent supporting errors therefore satisfy
$\sum_{i=0}^{P-1}\E[C_i+E_{i+1}]\le27B(P+2A)$. For all pairs $i<j$, the argument for
\eqref{eq:geom-triangle} bounds the expected sum of the negative parts of $\langle v_i,y_j\rangle$
by $\sum_{i=0}^P\E[C_i+E_i]\le54B(P+A)$. In the signed sum, the terms from $m_j$ contribute at
most $\sqrt B(P+2\sqrt A)$, since the weights $r^{i-j}$ sum to at most one over $i<j$.
The displacement terms contribute at most $\tfrac12\sum_{i=0}^P\E[\|v_i\|^2]\le4B(P+A)$.
Adding twice the negative parts to the signed sum in \eqref{eq:geom-triangle} gives
\[
 \E\!\left[\sum_{0\le i<j\le P}|\langle v_i,y_j\rangle|\right]\le128B(P+A).
\]

\subsubsection{Summing over scales and choosing the stepsize}

Let $\Phi\coloneqq\Phi_{n,2+\varepsilon}$. We first bound the rises of the envelope on the
dyadic grid in Appendix~\ref{app:scale-tails}. At time $n$, $I_{\lambda,n}\le a_\lambda$ for
every integer $\lambda\ge1$. The rise from $\lambda$ to $\lambda/2$ is at most
$R_\lambda(X_n)^2/(2\lambda)$. Squaring the three-term bound for $R_\lambda(X_n)/\sqrt\lambda$
and using $(u+v+w)^2\le3(u^2+v^2+w^2)$ gives a first part with second-moment norm at most
$3\cdot36^2/2<2^{11}$. The remainder has expectation at most
$3456\lambda^{-\varepsilon/2}+216A(\lambda/n)^{1/3}$. For $\lambda<1$ the whole rise is at most
$\lambda/2$. On the grid from $\Lambda\coloneqq2^{\lceil\log_2(n-1)\rceil}$ down to $1$,
the weights $\lambda^{-\varepsilon/2}$ sum to at most $3\Phi$ by the integral comparison in
Appendix~\ref{app:scale-tails}. The weights $(\lambda/n)^{1/3}$ sum to at most
$2^{1/3}/(1-2^{-1/3})<7$.

Each coarse rise contains $2\log_2 r$ grid increments. Thus the retained part $\delta_i'$ defined
in Appendix~\ref{app:scale-tails} satisfies $\|\delta_i'\|_2\le2^{13}\rho$. The initial range
contains at most $4+2\log_2 r$ increments, so the sum of its first parts also has norm at most
$2^{13}\rho$. Since $\E[M_\Lambda(X_n)]\le A+1$ and $f(X_n)-M_{1/8}(X_n)\le1/16$, summing gives
\begin{equation}\label{eq:numerical-grid}
 \E[f(X_n)]\le 2^{11}A+2^{14}\Phi+2^{13}\rho+\sum_{i=0}^{P-1}\E[\delta_i'] .
\end{equation}

For $0\le i<P$ let $B_i$ be the event that $\|v_i\|>r$ or $\|y_{i+1}\|>r$.
Applying the fourth- and second-moment bounds for the two vectors at threshold $r/3$ and using
$\mathbf1_{\{144a_N>r/3\}}\le432a_N/r$ gives
\[
 \sqrt{\Pp[B_i]}\le\frac{2^{21}}{r^2}
        +\frac{2^{11}N_{i+1}^{-\varepsilon/4}}r+\frac{2^{10}a_{N_i}}r .
\]
Cauchy-Schwarz bounds $\E[\delta_i'\mathbf1_{B_i}]$ by $2^{13}\rho$ times this expression.
The first terms sum to at most $P$ and the third to at most $\sqrt A$, since $r>2^{100}$ and
$\sum_{i=0}^P a_{N_i}\le2\sqrt A$. For the middle terms, the integral comparison in
Appendix~\ref{app:scale-tails} gives
$\rho\sum_{i=0}^{P-1}N_{i+1}^{-\varepsilon/4}\le\rho+2\Phi$.
As $2^{24}\rho/r\le1$ and $2^{25}/r\le1$, their sum costs at most $1+\Phi$.
The total contribution of the events $B_i$ is therefore at most $2(P+A+\Phi)$.

On $B_i^c$ we bound $\delta_i'$ by $C_i+E_{i+1}$ plus $2r$ times the absolute inner product of
the truncated vectors $\widetilde v_i,\widetilde y_{i+1}$ from Appendix~\ref{app:determinant}.
For the matrices $A_i$ defined there, expanding the quadratic forms and using the bound for all
pairs gives
\[
 \E\!\left[\sum_{i=0}^{P-1}\widetilde y_{i+1}^{\top}A_i\widetilde y_{i+1}\right]
 \le r^2\bigl(9B(P+2A)+128B(P+A)\bigr)\le160Br^2(P+A).
\]
The sum of the squared inverse-matrix norms of $\widetilde v_i$ is at most $2d\log(1+P/d)$.
Cauchy-Schwarz therefore bounds the expected contribution of the inner products, including
their coefficient $2r$, by $2\sqrt{320B}\,r^2\sqrt{(P+A)d\log(1+P/d)}\le2^{13}(P+A)$,
using \eqref{eq:numerical-scales}. Adding the adjacent supporting errors $27B(P+2A)$ and the
exceptional contributions in \eqref{eq:numerical-grid}, and using $\rho\le P$, gives
\[
 \E[f(X_n)]\le 2^{20}A+2^{19}P+2^{15}\Phi .
\]

Finally $P\le5q_{n,d}$, and $\Phi_{n,2+\varepsilon}\le2\Phi_{n,p}$ since $\Phi_{n,3}\le2$.
Restoring the stepsize as in Appendix~\ref{app:geometric} replaces $A$ by $D^2/(\eta n)$ in the
error bound and multiplies the noise terms by $\eta$. When $D\ge\eta n$, the elementary estimate
$D+2\eta(n-1)\le3D$ suffices. Combining these bounds with the scalar, small-moment, and
$w\le e^{1024}$ cases gives, for all parameters in Theorem~\ref{thm:main},
\begin{equation}\label{eq:numerical-upper}
 \mathcal R_{n,p}^{(d)}(\eta,D)
 \le 2^{20}a_\eta+(5\cdot2^{19}+2^{16})\eta b_{n,p,d}.
\end{equation}
This proves the upper bound in \eqref{eq:numerical-risk}.

For $D>0$ the argument in Appendix~\ref{app:main} gives
$a_\eta+\eta b_{n,p,d}\ge D\sqrt{b_{n,p,d}/n}$ for every $\eta>0$. At
$\eta=D/\sqrt{n b_{n,p,d}}$, each of the two terms equals $D\sqrt{b_{n,p,d}/n}$.
Since $2^{20}+5\cdot2^{19}+2^{16}=3735552<2^{22}$, \eqref{eq:numerical-upper} gives the same
upper comparison constant for this stepsize. Thus the lower and upper comparison constants
$1/43004$ and $2^{22}$ also hold in \eqref{eq:main-optimized}.

\end{document}